\documentclass[12pt,a4paper]{article}
\usepackage{amsfonts,amsmath,mathrsfs,amssymb, indentfirst}
\usepackage{tikz}
\usepackage{verbatim}
\usepackage{array}
\usepackage{setspace}
\newtheorem{theorem}{Theorem}[section]
\newtheorem{lemma}[theorem]{Lemma}

\newtheorem{defin}{Definition}[section]

\newtheorem{remark}{Remark}[section]

\begin{document}

\title {Symmetric versions of Laman's Theorem}

\author{Bernd Schulze\footnote{Research for this article was supported, in part, under a grant from NSERC (Canada), and final preparation occured at the TU Berlin with support of the DFG Research Unit 565 `Polyhedral Surfaces'.}\\ Inst. Mathematics, MA 6-2\\ TU Berlin\\ D-10623 Berlin, Germany
}

\maketitle

\begin{abstract}

Recent work has shown that if an isostatic bar and joint framework possesses non-trivial symmetries, then it must satisfy some very simply stated restrictions on the number of joints and bars that are `fixed' by various symmetry operations of the framework.\\\indent For the group $\mathcal{C}_3$ which describes 3-fold rotational symmetry in the plane, we verify the conjecture proposed in \cite{cfgsw} that these restrictions on the number of fixed structural components, together with the Laman conditions, are also  sufficient for a framework with $\mathcal{C}_3$ symmetry to be isostatic, provided that its joints are positioned as generically as possible subject to the given symmetry constraints.\\\indent In addition, we establish symmetric versions of Henneberg's Theorem and Crapo's Theorem for $\mathcal{C}_3$ which provide alternate characterizations of `generically' isostatic graphs with $\mathcal{C}_3$ symmetry.\\\indent As shown in \cite{BS4}, our techniques can be extended to establish analogous results for the symmetry groups $\mathcal{C}_2$ and $\mathcal{C}_s$ which are generated by a half-turn and a reflection in the plane, respectively.

\end{abstract}

\section{Introduction}
A bar and joint framework is said to be \emph{isostatic} if it is minimal infinitesimally rigid, in the sense that it is infinitesimally rigid and the removal of any bar results in an infinitesimally flexible framework (see \cite{graver, gss, W1, W2}, for example).\\\indent In 1970, G. Laman  provided an elegant characterization of \emph{generically $2$-isostatic} graphs, that is, graphs whose generic $2$-dimensional realizations as bar and joint frameworks are isostatic \cite{Lamanbib}. There are well known difficulties in extending this result to higher dimensions (see \cite{graver, gss, W2}, for example).\\\indent Using techniques from group representation theory, it was recently shown in \cite{cfgsw} that if a 2-dimensional isostatic  bar and joint framework possesses non-trivial symmetries, then it must not only satisfy the Laman conditions, but also some very simply stated extra conditions concerning the number of joints and bars that are fixed by various symmetry operations of the framework (see also \cite{BS2, BS4}). In particular, these restrictions imply that a $2$-dimensional isostatic framework must belong to one of only six possible point groups. In the Schoenflies notation \cite{bishop}, these groups are denoted by $\mathcal{C}_{1},\mathcal{C}_{2},\mathcal{C}_{3},\mathcal{C}_{s},\mathcal{C}_{2v}$, and $\mathcal{C}_{3v}$.\\\indent It was conjectured in \cite{cfgsw} that the Laman conditions, together with the corresponding additional conditions concerning the number of fixed structural components, are not only necessary, but also sufficient for a symmetric framework to be isostatic, provided that its joints are positioned as generically as possible subject to the given symmetry constraints.\\\indent
In this paper, we use the definition of `generic' for symmetry groups established in \cite{BS1} to verify this conjecture for the symmetry group $\mathcal{C}_3$ which describes 3-fold rotational symmetry in the plane ($\mathbb{Z}_3$ as an abstract group). The result is striking in its simplicity: to test a `generic' framework with $\mathcal{C}_3$ symmetry for isostaticity, we just need to check the number of joints that are `fixed' by the 3-fold rotation, as well as the standard conditions for generic rigidity without symmetry.\\\indent
By defining appropriate symmetrized inductive construction techniques, as well as appropriate symmetrized tree partitions of graphs, we also establish symmetric versions of Henneberg's Theorem (see \cite{gss, henne}) and Crapo's Theorem (\cite{Crapo, gss, tay}) for the group $\mathcal{C}_3$. These results provide us with some alternate techniques to give a `certificate' that a graph is `generically' isostatic modulo $\mathcal{C}_3$ symmetry. Furthermore, they enable us to generate all such graphs by means of an inductive construction sequence.
\\\indent With each of the main results presented in this paper, we also lay the foundation to design algorithms that decide whether a given graph is generically isostatic modulo $\mathcal{C}_3$ symmetry.\\\indent
It is shown in \cite{BS4} that our techniques can be extended to establish symmetric versions of Laman's Theorem, Henneberg's Theorem, and Crapo's Theorem for the symmetry groups $\mathcal{C}_{2}$ and $\mathcal{C}_{s}$ which are generated by a half-turn and a reflection, respectively, as well. However, it turns out that these proofs, in particular the ones for $\mathcal{C}_{s}$, are considerably more complex than the ones for $\mathcal{C}_{3}$. For simplicity, we therefore restrict our attention to the group $\mathcal{C}_{3}$ in this paper.\\\indent
The Laman-type conjectures for the dihedral groups $\mathcal{C}_{2v}$ and $\mathcal{C}_{3v}$ are still open. For a discussion on the difficulties that arise in proving these  conjectures (as well as a variety of related conjectures), we refer the interested reader to \cite{BS4}.

\section{Rigidity theoretic definitions and preliminaries}

\subsection{Graph theory terminology}

All graphs considered in this paper are finite graphs without loops or multiple edges. The \emph{vertex set} of a graph $G$ is denoted by $V(G)$ and the \emph{edge set} of $G$ is denoted by $E(G)$. Two vertices $u \ne v$ of $G$ are said to be \emph{adjacent} if $\{u,v\}\in E(G)$, and \emph{independent} otherwise. A set $S$ of vertices of $G$ is \emph{independent} if every two vertices of $S$ are independent. The \emph{neighborhood} $N_{G}(v)$ of a vertex $v\in V(G)$ is the set of all vertices that are adjacent to $v$ and the elements of $N_{G}(v)$ are called the \emph{neighbors} of $v$.

A graph $H$ is a \emph{subgraph} of $G$ if $V(H)\subseteq V(G)$ and $E(H)\subseteq E(G)$, in which case we write $H\subseteq G$.
The simplest type of subgraph of $G$ is that obtained by deleting a vertex or an edge from $G$. Let $v$ be a vertex and $e$ be an edge of $G$. Then we write $G-\{v$\} for the subgraph of $G$ that has $V(G)\setminus \{v\}$ as its vertex set and whose edges are those of $G$ that are not incident with $v$. Similarly, we write $G-\{e\}$ for the subgraph of $G$ that has $V(G)$ as its vertex set and $E(G)\setminus \{e\}$ as its edge set. The deletion of a set of vertices or a set of edges from $G$ is defined and denoted analogously.\\\indent
If $u$ and $v$ are independent vertices of $G$, then we write $G+\big\{\{u,v\}\big\}$ for the graph that has $V(G)$ as its vertex set and $E(G)\cup \big\{\{u,v\}\big\}$ as its edge set. The addition of a set of edges is again defined and denoted analogously.\\\indent
For a nonempty subset $U$ of $V(G)$, the subgraph $\langle U \rangle$ of $G$ \emph{induced} by $U$ is the graph having vertex set $U$ and whose edges are those of $G$ that are incident with two elements of $U$.\\\indent
The \emph{intersection} $G=G_{1}\cap G_{2}$ of two graphs $G_1$ and $G_2$ is the graph with $V(G)= V(G_{1})\cap V(G_{2})$ and $E(G)= E(G_{1})\cap E(G_{2})$. Similarly, the \emph{union} $G=G_{1}\cup G_{2}$ is the graph with $V(G)= V(G_{1})\cup V(G_{2})$ and $E(G)= E(G_{1})\cup E(G_{2})$.

An \emph{automorphism} of a graph $G$ is a permutation $\alpha$ of $V(G)$ such that $\{u,v\}\in E(G)$ if and only if $\{\alpha(u),\alpha(v)\}\in E(G)$.
The automorphisms of a graph $G$ form a group under composition which is denoted by $\textrm{Aut}(G)$.\\\indent
Let $H$ be a subgraph of $G$ and $\alpha\in \textrm{Aut}(G)$. We define $\alpha(H)$ to be the subgraph of $G$ that has $\alpha\big(V(H)\big)$ as its vertex set and $\alpha\big(E(H)\big)$ as its edge set, where $\{u,v\}\in \alpha\big(E(H)\big)$ if and only if $\alpha^{-1}(\{u,v\})=\{\alpha^{-1}(u),\alpha^{-1}(v)\}\in E(H)$.\\\indent We say that $H$ is \emph{invariant under} $\alpha$ if $\alpha\big(V(H)\big)=V(H)$ and $\alpha\big(E(H)\big)=E(H)$, in which case we write $\alpha(H)=H$.

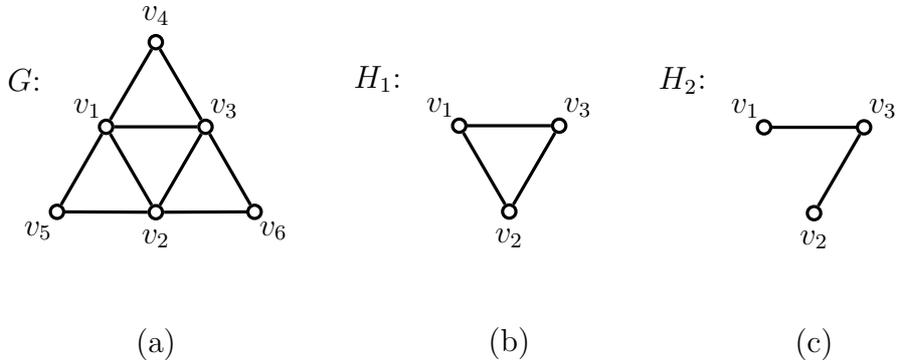
\begin{figure}[htp]
\begin{center}
\begin{tikzpicture}[very thick,scale=1]
\tikzstyle{every node}=[circle, draw=black, fill=white, inner sep=0pt, minimum width=5pt];
\path (30:0.755cm) node (v3) [label = above right: $v_3$] {} ;
\path (150:0.755cm) node (v1) [label = above left: $v_1$] {} ;
\path (270:0.755cm) node (v2) [label = below: $v_2$] {} ;
\path (90:1.5cm) node (v4) [label = above: $v_4$] {} ;
\path (210:1.5cm) node (v5) [label = below left: $v_5$] {} ;
\path (330:1.5cm) node (v6) [label = below right: $v_6$] {} ;
\draw (v1) -- (v2);
\draw (v1) -- (v3);
\draw (v3) -- (v2);
\draw (v1) -- (v4);
\draw (v4) -- (v3);
\draw (v1) -- (v5);
\draw (v2) -- (v5);
\draw (v6) -- (v2);
\draw (v6) -- (v3);
\node [draw=white, fill=white] (c) at (150:2cm) {$G$:};
\node [draw=white, fill=white] (c) at (270:2.5cm) {(a)};
 \end{tikzpicture}
 \hspace{0.5cm}
 \begin{tikzpicture}[very thick,scale=1]
\tikzstyle{every node}=[circle, draw=black, fill=white, inner sep=0pt, minimum width=5pt];
\path (30:0.759cm) node (v3) [label = above right: $v_3$] {} ;
\path (150:0.759cm) node (v1) [label = above left: $v_1$] {} ;
\path (270:0.759cm) node (v2) [label = below: $v_2$] {} ;
\draw (v1) -- (v2);
\draw (v1) -- (v3);
\draw (v3) -- (v2);
\node [draw=white, fill=white] (c) at (150:2cm) {$H_1$:};
\node [draw=white, fill=white] (c) at (270:2.5cm) {(b)};
\end{tikzpicture}
  \hspace{0.5cm}
 \begin{tikzpicture}[very thick,scale=1]
\tikzstyle{every node}=[circle, draw=black, fill=white, inner sep=0pt, minimum width=5pt];
\path (30:0.759cm) node (v3) [label = above right: $v_3$] {} ;
\path (150:0.759cm) node (v1) [label = above left: $v_1$] {} ;
\path (270:0.759cm) node (v2) [label = below: $v_2$] {} ;
\draw (v3) -- (v2);
\draw (v1) -- (v3);
\node [draw=white, fill=white] (c) at (150:2cm) {$H_2$:};
\node [draw=white, fill=white] (c) at (270:2.5cm) {(c)};
\end{tikzpicture}
         \end{center}
         \vspace{-0.3cm}
         \caption{\emph{ An invariant \emph{(b)} and a non-invariant subgraph \emph{(c)} of the graph $G$ under $\alpha=(v_{1}\, v_{2}\, v_{3})(v_{4}\, v_{5}\, v_{6})\in \textrm{Aut}(G)$.}}
\label{fig:invsubgraphs}
        \end{figure}

The graph $G$ in Figure \ref{fig:invsubgraphs} (a), for example, has the automorphism $\alpha=(v_{1}\, v_{2}\, v_{3})(v_{4}\, v_{5}\, v_{6}).$ The subgraph $H_{1}$ of $G$ is invariant under $\alpha$, but the subgraph $H_{2}$ of $G$ is not, because $\alpha\big(E(H_{2})\big)\neq E(H_{2})$.

\subsection{Infinitesimal rigidity}

\begin{defin}
\label{framework}
\emph{\cite{graver, gss, W1, W2} A \emph{framework} in $\mathbb{R}^{d}$ is a pair $(G,p)$, where $G$ is a graph and $p: V(G)\to \mathbb{R}^{d}$ is a map with the property that $p(u) \neq p(v)$ for all $\{u,v\} \in E(G)$. We also say that $(G,p)$
is a $d$-dimensional \emph{realization} of the \emph{underlying graph} $G$.\\\indent An ordered pair $\big(v,p(v)\big)$, where $v \in V(G)$, is a \emph{joint} of $(G,p)$, and an unordered pair $\big\{\big(u,p(u)\big),\big(v,p(v)\big)\big\}$ of joints, where $\{u,v\} \in E(G)$, is a \emph{bar} of $(G,p)$. }
\end{defin}

\begin{defin}
\label{infinmotion}
\emph{\cite{gss, W1, W2} Let $(G,p)$ be a framework in $\mathbb{R}^{d}$ with $V(G)=\{v_{1},v_{2},\ldots, v_{n}\}$. An \emph{infinitesimal motion} of $(G,p)$ is a function $u: V(G)\to \mathbb{R}^{d}$ such that
\begin{equation}
\label{infinmotioneq}
\big(p(v_{i})-p(v_{j})\big)\cdot \big(u(v_{i})-u(v_{j})\big)=0 \quad\textrm{ for all } \{v_{i},v_{j}\} \in E(G)\textrm{.}\nonumber\end{equation}
}
\end{defin}

\begin{figure}[htp]
\begin{center}
\begin{tikzpicture}[very thick,scale=1]
\tikzstyle{every node}=[circle, draw=black, fill=white, inner sep=0pt, minimum width=5pt];
        \path (0,0) node (p1) [label = below left: $p_{1}$] {} ;
        \path (2,0) node (p2) [label = below left: $p_2$] {} ;
        \draw (p1)  --  (p2);
        \draw [dashed, thin] (0.4,0) -- (0.4,1);
        \draw [dashed, thin] (2.4,0) -- (2.4,-0.4);
        \draw [dashed, thin] (2.4,0) -- (p2);
        \draw [->, black!60!white, ultra thick] (p1) -- node [draw=white, left=4pt] {$u_{1}$} (0.4,1);
        \draw [->, black!60!white, ultra thick] (p2) -- node [draw=white, below =4pt] {$u_2$} (2.4,-0.4);
        \node [draw=white, fill=white] (a) at (1,-2.5) {(a)};
        \end{tikzpicture}
        \hspace{0.5cm}
        \begin{tikzpicture}[very thick,scale=1]
\tikzstyle{every node}=[circle, draw=black, fill=white, inner sep=0pt, minimum width=5pt];
    \path (0,0) node (p1) [label = below left: $p_1$] {} ;
    \path (2,0) node (p2) [label = below right: $p_2$] {} ;
    \node (p3) at (1,0.1) {};
    \node [draw=white, fill=white] (labelp3) at (1,-0.25) {$p_{3}$};
    \draw (p1) -- (p2);
    \draw (0,0.1) -- (p3);
    \draw (p3) -- (2,0.1);
    \draw [->, black!60!white, ultra thick] (p3) -- node [draw=white, left=4pt] {$u_{3}$} (1,1);
    \draw [ black!60!white, ultra thick] (p1) -- node [rectangle, draw=white, above left] {$u_{1}=0$} (p1);
    \draw [ black!60!white, ultra thick] (p2) -- node [rectangle, draw=white, above right] {$u_{2}=0$} (p2);
    \node [draw=white, fill=white] (b) at (1,-2.5) {(b)};
    \end{tikzpicture}
    \hspace{0.5cm}
        \begin{tikzpicture}[very thick,scale=1]
\tikzstyle{every node}=[circle, draw=black, fill=white, inner sep=0pt, minimum width=5pt];
    \path (160:1.2cm) node (p6) [label = left: $p_6$] {} ;
    \path (120:1.2cm) node (p1) [label = left: $p_1$] {} ;
    \path (255:1.2cm) node (p2) [label = left: $p_2$] {} ;
     \path (20:1.2cm) node (p3) [label = below right: $p_3$] {} ;
    \path (60:1.2cm) node (p4) [label = right: $p_4$] {} ;
    \path (285:1.2cm) node (p5) [label = right: $p_5$] {} ;
     \draw (p1) -- (p4);
     \draw [->, black!60!white, ultra thick] (p6) --  (160:0.7)node[draw=white, below right=0.3pt]{$u_{6}$};
     \draw (p1) -- (p5);
     \draw (p1) -- (p6);
     \draw (p2) -- (p4);
     \draw (p2) -- (p5);
     \draw (p2) -- (p6);
     \draw (p3) -- (p4);
          \draw (p3) -- (p6);
     \draw [->, black!60!white, ultra thick] (p1) -- node[rectangle, draw=white, above right=3pt] {$u_{1}$}(120:1.7);
     \draw [->, black!60!white, ultra thick] (p2) -- node[rectangle, draw=white, below right=3pt] {$u_{2}$}(255:1.7);
     \draw [->, black!60!white, ultra thick] (p3) -- node[rectangle, draw=white, above=4pt] {$u_{3}$}(20:1.7);
     \draw [->, black!60!white, ultra thick] (p4)  -- node[rectangle, draw=white, left=4pt] {$u_{4}$} node[rectangle, draw=white, below=27pt] {$u_{5}$}(60:0.7);
     \draw (p3) -- (p5);
     \draw [->, black!60!white, ultra thick] (p5) -- (285:0.7);
     \node [draw=white, fill=white] (c) at (0,-2.5) {(c)};
        \end{tikzpicture}
\end{center}
\vspace{-0.3cm}
\caption{\emph{The arrows indicate the non-zero displacement vectors of an infinitesimal rigid motion \emph{(a)} and infinitesimal flexes \emph{(b, c)} of frameworks in $\mathbb{R}^2$.}}
\label{inmo}
\end{figure}
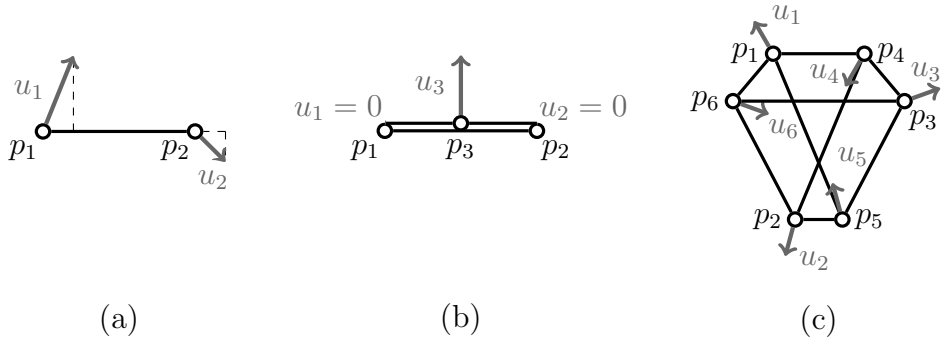

An infinitesimal motion is a set of displacement vectors, one at each joint of the framework, which preserve the lengths of all bars at first order (see also Figure \ref{inmo}).

\begin{defin}
\label{infinrigmotion}
\emph{\cite{gss, W1, W2} An infinitesimal motion $u$ of a framework $(G,p)$ is an \emph{infinitesimal rigid motion} if there exists a skew-symmetric matrix $S$ (a rotation) and a vector $t$ (a translation) such that $u(v)=Sp(v)+t$ for all $v\in V(G)$. Otherwise $u$ is an \emph{infinitesimal flex} of $(G,p)$.
}
\end{defin}

\begin{defin}
\emph{\cite{gss, W1, W2} A framework $(G,p)$ is \emph{infinitesimally rigid} if every infinitesimal motion of $(G,p)$ is an infinitesimal rigid motion. Otherwise $(G,p)$ is said to be \emph{infinitesimally flexible}.}
\end{defin}

For a framework $(G,p)$ whose underlying graph $G$ has a vertex set that is indexed from 1 to $n$, say $V(G)=\{v_{1},v_{2},\ldots ,v_{n}\}$, we will frequently denote $p(v_{i})$ by $p_{i}$ for $i=1,2,\ldots, n$. The $k^{th}$ component of a vector $x$ is denoted by $(x)_{k}$.

\begin{defin}
\emph{\cite{graver, gss, W1, W2} Let $G$ be a graph with $V(G)=\{v_{1},v_{2},\ldots,v_{n}\}$ and let $p:V(G)\to \mathbb{R}^{d}$. The
\emph{rigidity matrix} of $(G,p)$ is the $|E(G)| \times dn$ matrix  \begin{displaymath} \mathbf{R}(G,p)=\left(
\begin{array} {ccccccccccc }
& & & & & \vdots & & & & & \\
0 & \ldots & 0 & p_{i}-p_{j}&0 &\ldots &0 & p_{j}-p_{i} &0 &\ldots &
0\\ & & & & & \vdots & & & & &\end{array}
\right)\textrm{,}\end{displaymath} that is, for each edge $\{v_{i},v_{j}\}\in E(G)$, $\mathbf{R}(G,p)$ has the row with
$(p_{i}-p_{j})_{1},\ldots,(p_{i}-p_{j})_{d}$ in the columns $d(i-1)+1,\ldots,di$, $(p_{j}-p_{i})_{1},\ldots,(p_{j}-p_{i})_{d}$ in
the columns $d(j-1)+1,\ldots,dj$, and $0$ elsewhere.}
\end{defin}

Note that if we identify an infinitesimal motion of $(G,p)$ with a column vector in $\mathbb{R}^{dn}$ (by using the order on $V(G)$), then the kernel of the rigidity matrix $\mathbf{R}(G,p)$ is the space of all infinitesimal motions of $(G,p)$.

\begin{theorem}
\label{infinrigaff}\cite{asiroth, gluck}
A framework $(G,p)$ in $\mathbb{R}^d$ is infinitesimally rigid if and only if either $\textrm{rank }\big(\mathbf{R}(G,p)\big)=d |V(G)| - \binom{d+1}{2}$ or $G$ is a complete graph $K_n$ and the points $p(v)$, $v\in V(G)$, are affinely independent.
\end{theorem}

\begin{defin}
\emph{\cite{gss, W1, W2} A framework $(G,p)$ is \emph{independent} if the row vectors of the rigidity matrix $\mathbf{R}(G,p)$ are linearly independent.  A framework which is both independent and infinitesimally rigid is called \emph{isostatic}.
}
\end{defin}

An isostatic framework is minimal infinitesimally rigid, in the sense that the removal of any bar results in an infinitesimally flexible framework (see also \cite{graver, gss, W1, W2}).

\begin{theorem}
\label{isoequiv}\cite{gss, W2}
For a $d$-dimensional realization $(G,p)$ of a graph $G$ with $|V(G)|\geq d$, the following are equivalent:
\begin{itemize}
\item[(i)] $(G,p)$ is isostatic;
\item[(ii)] $(G,p)$ is infinitesimally rigid and $|E(G)|=d|V(G)|-\binom{d+1}{2}$;
\item[(iii)] $(G,p)$ is independent and $|E(G)|=d|V(G)|-\binom{d+1}{2}$;
\end{itemize}
\end{theorem}

\subsection{Generic rigidity}
\label{subsec:genrig}

Generic rigidity is concerned with the infinitesimal rigidity of `almost all' geometric realizations of a given graph. One of the `standard' definitions of `generic' that is frequently used in rigidity theory is the following.

\begin{defin}
\label{generic}
\emph{\cite{graver, gss} Let $G$ be a graph with $V(G)=\{v_1,\ldots, v_n\}$ and $K_n$ be the complete graph on $V(G)$. A framework $(G,p)$ is \emph{generic} if the determinant of any submatrix of $\mathbf{R}(K_n,p)$ is zero only if it is (identically) zero in the variables $p'_i$.}
\end{defin}

There are two fundamental facts regarding this definition of generic. First, it follows immediately from Definition \ref{generic} that the set of all generic realizations of a given graph $G$ in $\mathbb{R}^d$ forms a dense open subset of all possible realizations of $G$ in $\mathbb{R}^d$ \cite{gss}. Secondly, the infinitesimal rigidity properties are the same for all generic realizations of $G$, as the next result shows:

\begin{theorem}
\label{genericrigtheorem}\cite{graver, gss, W2}
For a graph $G$ and a fixed dimension $d$, the following are equivalent:
\begin{itemize}
\item[(i)] $(G,p)$ is infinitesimally rigid (independent, isostatic) for some map $p:V(G)\to \mathbb{R}^d$;
\item[(ii)] every $d$-dimensional generic realization of $G$ is infinitesimally rigid (independent, isostatic).
\end{itemize}
\end{theorem}
\noindent

It follows that for generic frameworks, infinitesimal rigidity is purely combinatorial, and hence a property of the underlying graph. This gives rise to the following definition of infinitesimal rigidity for graphs:

\begin{defin}
\emph{A graph $G$ is \emph{generically $d$-rigid ($d$-independent, $d$-isostatic)} if $d$-dimensional generic realizations of $G$ are infinitesimally rigid (independent, isostatic).}
\end{defin}

In 1970, G. Laman proved the following combinatorial characterization of generically $2$-isostatic graphs.

\begin{theorem}[Laman, 1970]
\label{lamantheorem}\cite{Lamanbib}
A graph $G$ with $|V(G)|\geq 2$ is generically 2-isostatic if and only if
\begin{itemize}
\item[(i)] $|E(G)|=2|V(G)|-3$;
\item[(ii)]$|E(H)|\leq 2|V(H)|-3$ for all $H\subseteq G$ with $|V(H)|\geq 2$.
\end{itemize}
\end{theorem}

Various proofs of Laman's Theorem can be found in \cite{graver}, \cite{gss}, \cite{lovyem}, \cite{tay}, and \cite{W7}, for example.\\\indent Throughout this paper, we will refer to the conditions $(i)$ and $(ii)$ in Theorem \ref{lamantheorem} as the \emph{Laman conditions}.

A combinatorial characterization of generically isostatic graphs in dimension 3 or higher is not yet known. See \cite{graver, gss, taywhit}, for example, for a detailed discussion of this problem.

There are some inductive construction techniques that preserve the generic rigidity properties of a graph. These construction techniques can be used to prove theorems such as Laman's Theorem, to analyze graphs for generic rigidity, and to characterize generically 1-isostatic and 2-isostatic graphs. For all dimensions $d$, they provide a tool to generate classes of generically $d$-isostatic graphs.

\begin{defin}
\label{ver2ex}
\emph{\cite{taywhit, W1} Let $G$ be a graph, $U\subseteq V(G)$ with $|U|=d$ and $v\notin V(G)$. Then the graph $\widehat{G}$ with $V(\widehat{G})=V(G)\cup \{v\}$ and $E(\widehat{G})=E(G)\cup\big\{\{v,u\}|u\in U\big\}$ is called a \emph{vertex $d$-addition (by $v$) of} $G$.}
\end{defin}

\begin{theorem}[Vertex Addition Theorem]\cite{graver, gss, taywhit, W1}
\label{vertext}
A vertex d-addition of a generically $d$-isostatic graph is generically $d$-isostatic. Conversely, deleting a vertex of valence $d$ from a generically $d$-isostatic graph results in a generically $d$-isostatic graph.
\end{theorem}

\begin{defin}
\emph{\cite{taywhit, W1} Let $G$ be a graph, $U\subseteq V(G)$ with $|U|=d+1$ and $\{u_{1},u_{2}\}\in E(G)$ for some $u_{1},u_{2}\in U$. Further, let $v\notin V(G)$. Then the graph $\widehat{G}$ with $V(\widehat{G})=V(G)\cup \{v\}$ and $E(\widehat{G})=\big(E(G)\setminus \big\{\{u_{1},u_{2}\}\big\}\big)\cup\big\{\{v,u\}|u\in U\big\}$ is called an \emph{edge $d$-split (on $u_{1},u_{2};v$) of $G$.}}
\end{defin}

\begin{theorem}[Edge Split Theorem]\cite{graver, gss, taywhit, W1}
\label{edgespl}
An edge $d$-split of a generically $d$-isostatic graph is generically $d$-isostatic. Conversely, if one deletes a vertex $v$ of valence $d+1$ from a generically $d$-isostatic graph, then one may add an edge between one of the pairs of vertices adjacent to $v$ so that the resulting graph is generically $d$-isostatic.
\end{theorem}

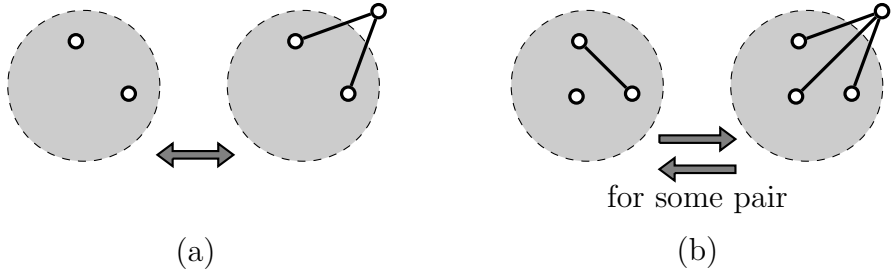
\begin{figure}[htp]
\begin{center}
\begin{tikzpicture}[very thick,scale=1]
\tikzstyle{every node}=[circle, draw=black, fill=white, inner sep=0pt, minimum width=5pt];
\filldraw[fill=black!20!white, draw=black, thin, dashed](0,0)circle(1cm);
\node (p1) at (100:0.6cm) {};
\node (p2) at (350:0.6cm) {};
  \end{tikzpicture}
\hspace{0.6cm}
\begin{tikzpicture}[very thick,scale=1]
\tikzstyle{every node}=[circle, draw=black, fill=white, inner sep=0pt, minimum width=5pt];
\filldraw[fill=black!20!white, draw=black, thin, dashed](0,0)circle(1cm);
\node (p1) at (100:0.6cm) {};
\node (p2) at (350:0.6cm) {};
\node (p3) at (45:1.4cm) {};
\draw (p3)--(p1);
\draw (p3)--(p2);
\end{tikzpicture}
\hspace{1.35cm}
\begin{tikzpicture}[very thick,scale=1]
\tikzstyle{every node}=[circle, draw=black, fill=white, inner sep=0pt, minimum width=5pt];
\filldraw[fill=black!20!white, draw=black, thin, dashed](0,0)circle(1cm);
\node (p1) at (100:0.6cm) {};
\node (p2) at (350:0.6cm) {};
\node (p4) at (225:0.2cm) {};
 \draw(p1)--(p2);
\end{tikzpicture}
\hspace{0.6cm}
\begin{tikzpicture}[very thick,scale=1]
\tikzstyle{every node}=[circle, draw=black, fill=white, inner sep=0pt, minimum width=5pt];
\filldraw[fill=black!20!white, draw=black, thin, dashed](0,0)circle(1cm);
\node (p1) at (100:0.6cm) {};
\node (p2) at (350:0.6cm) {};
\node (p3) at (45:1.4cm) {};
\node (p4) at (225:0.2cm) {};
\draw (p3)--(p1);
\draw (p3)--(p2);
\draw (p3)--(p4);
\end{tikzpicture}
\vspace{-0.5cm}

\begin{tikzpicture}[very thick,scale=1]
\tikzstyle{every node}=[circle, draw=black, fill=white, inner sep=0pt, minimum width=5pt];
\node [draw=white, fill=white] (a) at (1.9,-2.3) {(a)};
\node [draw=white, fill=white] (a) at (1,-1.2) {};
\filldraw[fill=black!50!white, draw=black, thick]
    (1.6,-1.05) -- (2.2,-1.05) -- (2.2,-1.15) -- (2.4,-1.0) -- (2.2,-0.85) -- (2.2,-0.95) -- (1.6,-0.95) -- (1.6,-0.85)
    -- (1.4,-1)--(1.6,-1.15)-- (1.6,-1.05);
\end{tikzpicture}
\hspace{4.6cm}
\begin{tikzpicture}[very thick,scale=1]
\tikzstyle{every node}=[circle, draw=black, fill=white, inner sep=0pt, minimum width=5pt];
\node [draw=white, fill=white] (b) at (1.5,-2.3) {(b)};
\filldraw[fill=black!50!white, draw=black, thick]
    (1,-0.85) -- (1.8,-0.85) -- (1.8,-0.95) -- (2,-0.8) -- (1.8,-0.65) -- (1.8,-0.75) -- (1,-0.75)-- (1,-0.85);
\filldraw[fill=black!50!white, draw=black, thick]
    (1.2,-1.25) --  (2,-1.25) -- (2,-1.15) -- (1.2,-1.15) -- (1.2,-1.05)
    -- (1,-1.2)--(1.2,-1.35)-- (1.2,-1.25);
\node [rectangle, draw=white, fill=white] (c) at (1.5,-1.6) {for some pair};
\end{tikzpicture}
\end{center}
\vspace{-0.3cm}
\caption{\emph{Illustrations of the Vertex Addition Theorem \emph{(a)} and the Edge Split Theorem \emph{(b)} in dimension 2.}}
\label{fig:inductive}
\end{figure}

In 1911, L. Henneberg gave the following characterization of generically 2-isostatic graphs.

\begin{theorem}[Henneberg, 1911]\cite{henne}
\label{henneberg}
A graph is generically 2-isostatic if and only if it may be constructed from a single edge by a sequence of vertex 2-additions and edge 2-splits.
\end{theorem}

For a proof of Henneberg's Theorem, see \cite{gss} or \cite{taywhit}, for example.

Another way of characterizing generically 2-isostatic graphs is due to H. Crapo and uses partitions of a graph into edge disjoint trees.

\begin{defin}
\label{3T2}
\emph{\cite{Crapo, LMW, tay} A \emph{3Tree2 partition} of a graph $G$ is a partition of $E(G)$ into the edge sets of three edge disjoint trees $T_{0},T_{1},T_{2}$ such that each vertex of G belongs to exactly two of the trees.\\\indent A 3Tree2 partition is called \emph{proper} if no non-trivial subtrees of distinct trees $T_{i}$ have the same span (i.e., the same vertex sets).}
\end{defin}

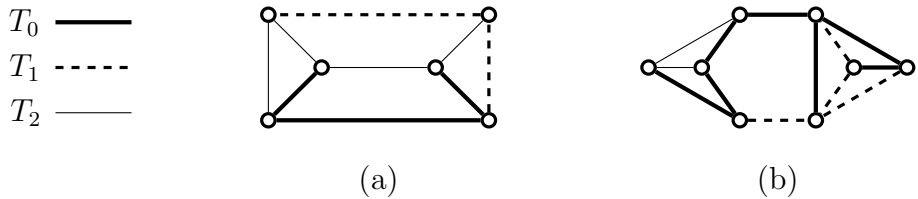
\begin{figure}[htp]
\begin{center}
\begin{tikzpicture}[very thick,scale=1]
\tikzstyle{every node}=[circle, draw=black, fill=white, inner sep=0pt, minimum width=5pt];
\draw [ultra thick] (-3.5,0.6) node[rectangle, draw=white,fill=white,left =5pt] {$T_0$} -- (-2.5,0.6);
\draw [dashed] (-3.5,0)  node[rectangle, draw=white,fill=white,left =5pt] {$T_1$} -- (-2.5,0);
\draw [thin] (-3.5,-0.6)  node[rectangle, draw=white,fill=white,left =5pt] {$T_2$} -- (-2.5,-0.6);
\node (p1) at (0,0) {};
\node (p2) at (1.5,0) {};
\node (p3) at (-0.7,-0.7) {};
\node (p4) at (-0.7,0.7) {};
\node (p5) at (2.2,-0.7) {};
\node (p6) at (2.2,0.7) {};
\draw[thin](p1)--(p2);
\draw [ultra thick] (p1)--(p3);
\draw[thin](p1)--(p4);
\draw [ultra thick](p2)--(p5);
\draw[thin](p2)--(p6);
\draw [ultra thick](p5)--(p3);
\draw[dashed](p6)--(p4);
\draw[dashed](p5)--(p6);
\draw[thin](p3)--(p4);
\node [draw=white, fill=white] (a) at (0.75,-1.5) {(a)};
\node (a1) at (5,0) {};
\node (a2) at (7,0) {};
\node (a3) at (4.3,0) {};
\node (a4) at (7.7,0) {};
\node (a5) at (5.5,0.7) {};
\node (a6) at (5.5,-0.7) {};
\node (a7) at (6.5,0.7) {};
\node (a8) at (6.5,-0.7) {};
\draw [ultra thick](a1)--(a5);
\draw [ultra thick](a1)--(a6);
\draw[thin](a1)--(a3);
\draw[thin](a3)--(a5);
\draw [ultra thick](a3)--(a6);
\draw [ultra thick](a7)--(a5);
\draw[dashed](a6)--(a8);
\draw [ultra thick](a7)--(a8);
\draw [ultra thick](a2)--(a4);
\draw[dashed](a2)--(a7);
\draw[dashed](a2)--(a8);
\draw [ultra thick](a4)--(a7);
\draw[dashed](a4)--(a8);
\node [draw=white, fill=white] (b) at (6,-1.5) {(b)};
\end{tikzpicture}
\end{center}
\vspace{-0.3cm}
\caption{\emph{A proper \emph{(a)} and a non-proper \emph{(b)} 3Tree2 partition.}}
\label{3T2fig}
\end{figure}

\begin{remark}
\label{3t2remk}
\emph{If a graph $G$ has a 3Tree2 partition, then it satisfies $|E(G)|= 2|V(G)|-3$. This follows from the presence of exactly two trees at each vertex of $G$ and the fact that for every tree $T$ we have $|E(T)|= |V(T)|-1$. Moreover, note that a 3Tree2 partition of a graph $G$ is proper if and only if every non-trivial subgraph $H$ of $G$ satisfies the count $|E(H)|\leq 2|V(H)|-3$ \cite{LMW}.}
\end{remark}

\begin{theorem}[Crapo, 1989]\cite{Crapo}
\label{3T2theorem}
A graph $G$ is generically 2-isostatic if and only if $G$ has a proper 3Tree2 partition.
\end{theorem}

\section{Symmetric frameworks}

A \emph{symmetry operation} of a framework $(G,p)$ in $\mathbb{R}^{d}$ is an isometry $x$ of $\mathbb{R}^{d}$ such that for some $\alpha\in \textrm{Aut}(G)$, we have
$x\big(p(v)\big)=p\big(\alpha(v)\big)$ for all $v\in V(G)$ \cite{Hall, BS1, BS2, BS4}.\\\indent
The set of all symmetry operations of a framework $(G,p)$ forms a group under composition, called the \emph{point group} of $(G,p)$ \cite{altherz, bishop, Hall, BS4, BS1}. Since translating a framework does not change its rigidity properties, we may assume wlog that the point group of any framework in this paper is a \emph{symmetry group}, i.e., a subgroup of the orthogonal group $O(\mathbb{R}^{d})$ \cite{BS2, BS4, BS1}.

We use the Schoenflies notation for the symmetry operations and symmetry groups considered in this paper, as this is one of the standard notations in the literature about symmetric structures (see \cite{altherz, bishop, cfgsw, FGsymmax, gsw, Hall, BS1, BS2, BS4}, for example). The three kinds of possible symmetry operations in dimension $2$ are the identity $Id$, rotations $C_{m}$ about the origin by an angle of $\frac{2\pi}{m}$, where $m\geq 2$, and reflections $s$ in lines through the origin. In the Schoenflies notation, this gives rise to the following families of possible symmetry groups in dimension 2: $\mathcal{C}_{1}$, $\mathcal{C}_{s}$, $\mathcal{C}_{m}$ and $\mathcal{C}_{mv}$, where $m\geq 2$. $\mathcal{C}_{1}$ denotes the trivial group which only contains the identity $Id$. $\mathcal{C}_{s}$ denotes any symmetry group in dimension 2 that consists of the identity $Id$ and a single reflection $s$. For $m\geq 2$, $\mathcal{C}_{m}$ denotes any cyclic symmetry group of order $m$ which is generated by a rotation $C_{m}$, and $\mathcal{C}_{mv} $ denotes any symmetry group in dimension 2 that is generated by a pair $\{C_{m},s\}$.

Given a symmetry group $S$ in dimension $d$ and a graph $G$, we let $\mathscr{R}_{(G,S)}$ denote the set of all $d$-dimensional realizations of $G$ whose point group is either equal to $S$ or contains $S$ as a subgroup \cite{BS1, BS2, BS4}. In other words, the set $\mathscr{R}_{(G,S)}$ consists of all realizations $(G,p)$ of $G$ for which there exists a map $\Phi:S\to \textrm{Aut}(G)$ so that
\begin{equation}\label{class} x\big(p(v)\big)=p\big(\Phi(x)(v)\big)\textrm{ for all } v\in V(G)\textrm{ and all } x\in S\textrm{.}\end{equation}

A framework $(G,p)\in \mathscr{R}_{(G,S)}$ satisfying the equations in (\ref{class}) for the map $\Phi:S\to \textrm{Aut}(G)$ is said to be \emph{of type $\Phi$}, and the set of all realizations in $\mathscr{R}_{(G,S)}$ which are of type $\Phi$ is denoted by $\mathscr{R}_{(G,S,\Phi)}$ (see again \cite{BS1, BS2, BS4}).

\begin{figure}[htp]
\begin{center}
\begin{tikzpicture}[very thick,scale=1]
\tikzstyle{every node}=[circle, draw=black, fill=white, inner sep=0pt, minimum width=5pt];
    \path (0.3,-0.5) node (p5) [label = below left: $p_{5}$] {} ;
    \path (1.5,-1) node (p3) [label = below left: $p_{3}$] {} ;
    \path (2.7,-0.5) node (p6) [label = below right: $p_{6}$] {} ;
   \path (0.8,0.7) node (p1) [label = left: $p_{1}$] {} ;
   \path (2.2,0.7) node (p2) [label = right: $p_{2}$] {} ;
   \path (1.5,1.1) node (p4) [label = above left: $p_{4}$] {} ;
   \draw (p1) -- (p4);
     \draw (p1) -- (p5);
     \draw (p1) -- (p6);
     \draw (p2) -- (p4);
     \draw (p2) -- (p5);
     \draw (p2) -- (p6);
     \draw (p3) -- (p4);
     \draw (p3) -- (p5);
     \draw (p3) -- (p6);
   \draw [dashed, thin] (1.5,-2) -- (1.5,2);
      \node [draw=white, fill=white] (a) at (1.5,-2.5) {(a)};
    \end{tikzpicture}
    \hspace{2cm}
        \begin{tikzpicture}[very thick,scale=1]
\tikzstyle{every node}=[circle, draw=black, fill=white, inner sep=0pt, minimum width=5pt];
    \path (160:1.2cm) node (p6) [label = left: $p_6$] {} ;
    \path (120:1.2cm) node (p1) [label = above left: $p_1$] {} ;
    \path (255:1.2cm) node (p2) [label = below left: $p_2$] {} ;
     \path (20:1.2cm) node (p3) [label = right: $p_3$] {} ;
    \path (60:1.2cm) node (p4) [label = above right: $p_4$] {} ;
    \path (285:1.2cm) node (p5) [label = below right: $p_5$] {} ;
     \draw (p1) -- (p4);
     \draw (p1) -- (p5);
     \draw (p1) -- (p6);
     \draw (p2) -- (p4);
     \draw (p2) -- (p5);
     \draw (p2) -- (p6);
     \draw (p3) -- (p4);
     \draw (p3) -- (p5);
     \draw (p3) -- (p6);
     \draw [dashed, thin] (0,-2) -- (0,2);
      \node [draw=white, fill=white] (b) at (0,-2.5) {(b)};
        \end{tikzpicture}
\end{center}
\vspace{-0.3cm}
\caption{\emph{$2$-dimensional realizations in $\mathscr{R}_{(K_{3,3},\mathcal{C}_s)}$ of different types.}}
\label{K33types}
\end{figure}
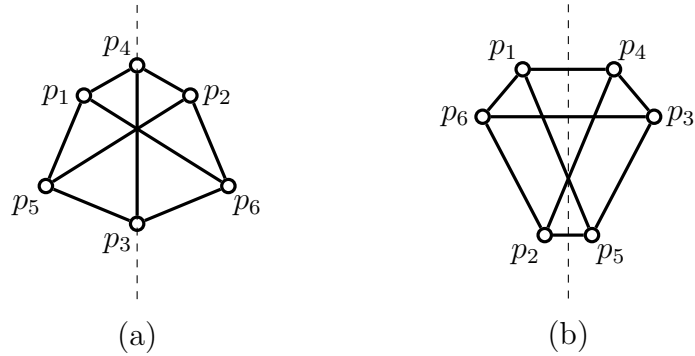

Different choices of types $\Phi:S\to \textrm{Aut}(G)$ frequently lead to very different geometric types of realizations of $G$ within $\mathscr{R}_{(G,S)}$. This is illustrated by the realizations of the complete bipartite graph $K_{3,3}$ with mirror symmetry depicted in Figure \ref{K33types}. The framework in Figure \ref{K33types} (a) is a realization in $\mathscr{R}_{(K_{3,3},\mathcal{C}_s)}$  of type
$\Phi_{a}$, where $\Phi_{a}: \mathcal{C}_{s} \to \textrm{Aut}(K_{3,3})$ is defined by
\begin{eqnarray} \Phi_{a}(Id)& =& id\nonumber\\\Phi_{a}(s)&=&
(v_{1}\,v_{2})(v_{5}\,v_{6})(v_{3})(v_{4})\textrm{,}\nonumber
\end{eqnarray}
and the framework in Figure \ref{K33types} (b)
is a realization in $\mathscr{R}_{(K_{3,3},\mathcal{C}_s)}$ of type
$\Phi_{b}$, where $\Phi_{b}: \mathcal{C}_{s} \to \textrm{Aut}(K_{3,3})$ is defined by
\begin{eqnarray}\Phi_{b}(Id)& = &id\nonumber\\ \Phi_{b}(s) &=& (v_{1}\,v_{4})(v_{2}\,v_{5})(v_{3}\,v_{6})\textrm{.}\nonumber
\end{eqnarray}
Note that `almost all' realizations in $\mathscr{R}_{(K_{3,3},\mathcal{C}_s,\Phi_a)}$  are isostatic, whereas all realizations in $\mathscr{R}_{(K_{3,3},\mathcal{C}_s,\Phi_b)}$ are infinitesimally flexible since the joints of any realization in $\mathscr{R}_{(K_{3,3},\mathcal{C}_s,\Phi_b)}$ are forced to lie on a conic section \cite{BS1, W3}.

\begin{remark}
\emph{A set $\mathscr{R}_{(G,S)}$ can possibly be empty. For example, there clearly exists no $2$-dimensional realization of $K_{2}$ in the set $\mathscr{R}_{(K_2,\mathcal{C}_{3})}$.\\\indent Given a non-empty set $\mathscr{R}_{(G,S)}$, it is also possible that $\mathscr{R}_{(G,S,\Phi)}=\emptyset$ for some map $\Phi:S\to \textrm{Aut}(G)$.\\\indent Consider, for example, the non-empty set $\mathscr{R}_{(K_{2},\mathcal{C}_{2})}$, where $\mathcal{C}_{2}=\{Id,C_{2}\}$ is the half-turn symmetry group in dimension 2, and let $I:\mathcal{C}_{2}\to \textrm{Aut}(K_{2})$ be the map which sends both $Id$ and $C_{2}$ to the identity automorphism of $K_{2}$. If $(K_{2},p)\in\mathscr{R}_{(K_{2},\mathcal{C}_{2},I)}$, then both joints of $(K_{2},p)$ must be located at the origin (which is the center of $C_{2}$). This contradicts Definition \ref{framework} of a framework, and hence we have $\mathscr{R}_{(K_{2},\mathcal{C}_{2},I)}=\emptyset$.}
\end{remark}

The following symmetry-adapted notion of generic for the set $\mathscr{R}_{(G,S,\Phi)}$ was introduced in \cite{BS1} (see also \cite{BS4}).

\begin{defin}
\label{symgeneric}
\emph{Let $G$ be a graph with $V(G)=\{v_1,\ldots, v_n\}$ and let $K_n$ be the complete graph on $V(G)$. Further, let $S$ be a symmetry group  and $\Phi$ be a map from $S$ to $\textrm{Aut}(G)$. A framework $(G,p)\in \mathscr{R}_{(G,S,\Phi)}$ is \emph{$(S,\Phi)$-generic} if the determinant of any submatrix of $\mathbf{R}(K_n,p)$ is zero only if it is zero for all $p'$ satisfying the symmetry equations in (\ref{class}).}
\end{defin}

Intuitively, an $(S,\Phi)$-generic realization of a graph $G$ is obtained by placing the vertices of a set of representatives for the symmetry orbits $S v=\{\Phi(x)(v)|\, x\in S\}$ into `generic' positions. The positions for the remaining vertices of $G$ are then uniquely determined by the symmetry constraints imposed by $S$ and $\Phi$ (see \cite{BS1, BS4}, for further details).

It is shown in \cite{BS1} that the set of $(S,\Phi)$-generic realizations of a graph $G$ is an open dense subset of the set $\mathscr{R}_{(G,S,\Phi)}$. Moreover, the infinitesimal rigidity properties are the same for all $(S,\Phi)$-generic realizations of $G$, as the following theorem shows.

\begin{theorem}
\label{symgenrigthm} \cite{BS1, BS4}
Let $G$ be a graph, $S$ be a symmetry group, and $\Phi$ be a map from $S$ to $\textrm{Aut}(G)$ such that $\mathscr{R}_{(G,S,\Phi)}\neq \emptyset$. The following are equivalent.
\begin{itemize}
\item[(i)] There exists a framework $(G,p)\in \mathscr{R}_{(G,S,\Phi)}$ that is infinitesimally rigid (independent, isostatic);
\item[(ii)] every $(S,\Phi)$-generic realization of $G$ is infinitesimally rigid (independent, isostatic).
\end{itemize}
\end{theorem}

So, being infinitesimally rigid (independent, isostatic) is an $(S,\Phi)$-generic property. This gives rise to

\begin{defin}
\emph{Let $G$ be a graph, $S$ be a symmetry group, and $\Phi$ be a map from $S$ to $\textrm{Aut}(G)$. Then $G$ is said to be \emph{$(S,\Phi)$-generically infinitesimally rigid (independent, isostatic)} if all realizations of $G$ which are $(S,\Phi)$-generic are infinitesimally rigid (independent, isostatic).}
\end{defin}

Using techniques from group representation theory, it is shown in \cite{cfgsw} that if a symmetric isostatic framework $(G,p)$ belongs to a set $\mathscr{R}_{(G,S,\Phi)}$, where $S$ is a non-trivial symmetry group and $\Phi:S\to \textrm{Aut}(G)$ is a homomorphism, then $(G,p)$ needs to satisfy certain restrictions on the number of joints and bars that are `fixed' by various symmetry operations of $(G,p)$ (see also \cite{FGsymmax,BS1, BS4}). In the following, we summarize the key result for dimension 2.

\begin{defin}
\emph{\cite{BS2, BS4} Let $G$ be a graph with $V(G)=\{v_{1},\ldots,v_{n}\}$, $S$ be a symmetry group, $\Phi$ be a map from $S$ to $\textrm{Aut}(G)$, $(G,p)$ be a framework in $\mathscr{R}_{(G,S,\Phi)}$, and $x\in S$. A joint $(v_{i},p_{i})$ of $(G,p)$ is said to be \emph{fixed by $x$ with respect to $\Phi$} if $\Phi(x)(v_{i})=v_{i}$.\\\indent Similarly, a bar $\{(v_{i},p_{i}),(v_{j},p_{j})\}$ of $(G,p)$ is said to be \emph{fixed by $x$ with respect to $\Phi$} if $\Phi(x)\big(\{v_{i},v_{j}\}\big)=\{v_{i},v_{j}\}$.\\\indent The number of joints of $(G,p)$ that are fixed by $x$ with respect to $\Phi$ is denoted by $j_{\Phi(x)}$ and the number of bars of $(G,p)$ that are fixed by $x$ with respect to $\Phi$ is denoted by $b_{\Phi(x)}$.}
\end{defin}

\begin{remark}
\emph{If a joint $\big(v,p(v)\big)$ of a framework $(G,p)\in\mathscr{R}_{(G,S,\Phi)}$ is fixed by $x\in S$ with respect to $\Phi$, then we have $x\big(p(v)\big)=p\big(\Phi(x)(v)\big)=p(v)$. In particular, if $(G,p)$ is a 2-dimensional framework, then a joint that is fixed by a rotation $C_m\in S$ must lie at the center of $C_m$, and a joint that is fixed by a reflection $s\in S$ must lie on the mirror line corresponding to $s$. Similar geometric restrictions of course also apply for \emph{bars} of $(G,p)$ that are fixed by various symmetry operations in $S$ (see \cite{cfgsw, BS1, BS4} for details).}
\end{remark}

\begin{theorem}
\label{symmaxrule} \cite{cfgsw, BS4}
Let $G$ be a graph, $S$ be a symmetry group in dimension $2$, $\Phi:S\to \textrm{Aut}(G)$ be a homomorphism, and $(G,p)$ be an isostatic framework in $\mathscr{R}_{(G,S,\Phi)}$ with the property that the points $p(v)$, $v\in V(G)$, span all of $\mathbb{R}^{2}$. Then
\begin{itemize}
\item[(i)] the Laman conditions are satisfied;
\item[(ii)] if $S=\mathcal{C}_{2}$, then  $j_{\Phi(C_{2})}=0$ and $b_{\Phi(C_{2})}=1$;
\item[(iii)] if $S=\mathcal{C}_{3}$, then  $j_{\Phi(C_{3})}=0$;
\item[(iv)] if $S=\mathcal{C}_{s}$, then  $b_{\Phi(s)}=1$;
\item[(v)]if $S=\mathcal{C}_{2v}$, then $j_{\Phi(C_{2})}=0$ and $b_{\Phi(C_{2})}=b_{\Phi(s)}=1$ for both reflections $s\in \mathcal{C}_{2v}$;
\item[(vi)]if $S=\mathcal{C}_{3v}$, then $j_{\Phi(C_{3})}=0$ and $b_{\Phi(s)}=1$ for all reflections $s\in \mathcal{C}_{3v}$;
\item[(vii)] $S$ is either the trivial group $\mathcal{C}_{1}$ or one of the five non-trivial symmetry groups listed above.
\end{itemize}
\end{theorem}

Examples of isostatic frameworks for each of the point groups listed in Theorem \ref{symmaxrule} are given in \cite{cfgsw, BS4}.

It was conjectured in \cite{cfgsw} that the conditions identified in Theorem \ref{symmaxrule} $(ii)$ - $(vi)$, together with the Laman conditions, are also \emph{sufficient} for an $(S,\Phi)$-generic realization of $G$ to be isostatic. In the following, we verify this conjecture for the symmetry group $\mathcal{C}_3$.  In addition, we provide Henneberg-type and Crapo-type characterizations of $(\mathcal{C}_3,\Phi)$-generically isostatic graphs. The techniques used in these proofs are extended in \cite{BS4}  to prove the corresponding Laman-type conjectures for the groups $\mathcal{C}_2$ and $\mathcal{C}_s$ as well as  analogous Henneberg-type and Crapo-type results for these groups. Characterizations of $(\mathcal{C}_{2v},\Phi)$- or $(\mathcal{C}_{3v},\Phi)$-generically isostatic graphs, however, have not yet been established (see again \cite{BS4}).\\\indent While, initially, the fact that $\mathcal{C}_{3}$  allows the easiest and most natural proof for the Laman-type conjecture (as well as for a symmetric version of Crapo's Theorem) came as somewhat of a surprise, we can now identify some clear indications for this.\\\indent For example, Crapo's Theorem uses partitions of the edges of $G$ into \emph{three} edge-disjoint trees, so that it is most natural to extend this result to the cyclic group $\mathcal{C}_{3}$ of order three. Moreover, the condition $j_{\Phi(C_3)}=0$ implies that for any subgraph $H$ of $G$ with full $\mathcal{C}_{3}$ symmetry we must have that both $|V(H)|$ and $|E(H)|$ are multiples of three, so that $H$ cannot satisfy the count $|E(H)|=2|V(H)|-4$ or $|E(H)|=2|V(H)|-5$. As we will see, this turns out to be extremely useful in the proof of the Laman-type result for $\mathcal{C}_{3}$.

\section{Symmetric Henneberg moves and 3Tree2 partitions for $\mathcal{C}_{3}$}
\label{sec:henne3t2}

We need the following inductive construction techniques to obtain a symmetrized Henneberg's Theorem for $\mathcal{C}_{3}$.

\begin{figure}[htp]
\begin{center}
\begin{tikzpicture}[very thick,scale=1]
\tikzstyle{every node}=[circle, draw=black, fill=white, inner sep=0pt, minimum width=5pt];
\filldraw[fill=black!20!white, draw=black, thin, dashed](0,0)circle(1.6cm);
\node[rectangle,draw=black!20!white,fill=black!20!white](l5) at (290:0.81cm){$\gamma(v_1)$};
\node[rectangle,draw=black!20!white,fill=black!20!white](l6) at (360:0.71cm){$\gamma(v_2)$};
\node (p1) at (70:1.3cm) {};
\node (p2) at (110:1.3cm) {};
\node (p3) at (190:1.3cm) {};
\node (p4) at (230:1.3cm) {};
\node (p5) at (310:1.3cm) {};
\node (p6) at (350:1.3cm) {};
\node[rectangle,draw=black!20!white,fill=black!20!white](l1) at (50:1.05cm){$\gamma^2(v_1)$};
\node[rectangle,draw=black!20!white,fill=black!20!white](l2) at (130:1.05cm){$\gamma^2(v_2)$};
\node[rectangle,draw=black!20!white,fill=black!20!white](l3) at (180:0.95cm){$v_1$};
\node[rectangle,draw=black!20!white,fill=black!20!white](l4) at (240:0.95cm){$v_2$};
\filldraw[fill=black!50!white, draw=black, thick]
    (2.75,0) -- (3.35,0) -- (3.35,-0.1) -- (3.55,0.05) -- (3.35,0.2) -- (3.35,0.1) -- (2.75,0.1) -- (2.75,0);
\end{tikzpicture}
\hspace{0.5cm}
\begin{tikzpicture}[very thick,scale=1]
\tikzstyle{every node}=[circle, draw=black, fill=white, inner sep=0pt, minimum width=5pt];
\filldraw[fill=black!20!white, draw=black, thin, dashed](0,0)circle(1.6cm);
\node[rectangle,draw=black!20!white,fill=black!20!white](l5) at (290:0.81cm){$\gamma(v_1)$};
\node[rectangle,draw=black!20!white,fill=black!20!white](l6) at (360:0.71cm){$\gamma(v_2)$};
\node (p1) at (70:1.3cm) {};
\node (p2) at (110:1.3cm) {};
\node (p3) at (190:1.3cm) {};
\node (p4) at (230:1.3cm) {};
\node (p5) at (310:1.3cm) {};
\node (p6) at (350:1.3cm) {};
\node (p7) at (90:2.2cm) {};
\node (p8) at (210:2.2cm) {};
\node (p9) at (330:2.2cm) {};
\node[rectangle,draw=black!20!white,fill=black!20!white](l1) at (50:1.05cm){$\gamma^2(v_1)$};
\node[rectangle,draw=black!20!white,fill=black!20!white](l2) at (130:1.05cm){$\gamma^2(v_2)$};
\node[rectangle,draw=black!20!white,fill=black!20!white](l3) at (180:0.95cm){$v_1$};
\node[rectangle,draw=black!20!white,fill=black!20!white](l4) at (240:0.95cm){$v_2$};
\node[rectangle,draw=white,fill=white](l7) at (100:2.3cm){$z$};
\node[rectangle,draw=white,fill=white](l8) at (210:2.52cm){$v$};
\node[rectangle,draw=white,fill=white](l9) at (330:2.52cm){$w$};
\draw(p7)--(p1);
\draw(p7)--(p2);
\draw(p8)--(p3);
\draw(p8)--(p4);
\draw(p9)--(p5);
\draw(p9)--(p6);
\end{tikzpicture}
\end{center}
\vspace{-0.3cm}
\caption{\emph{A $(\mathcal{C}_{3},\Phi)$ vertex addition of a graph $G$, where $\Phi(C_{3})=\gamma$ and $\Phi(C_{3}^2)=\gamma^2$.}}
\label{C3vex}
\end{figure}
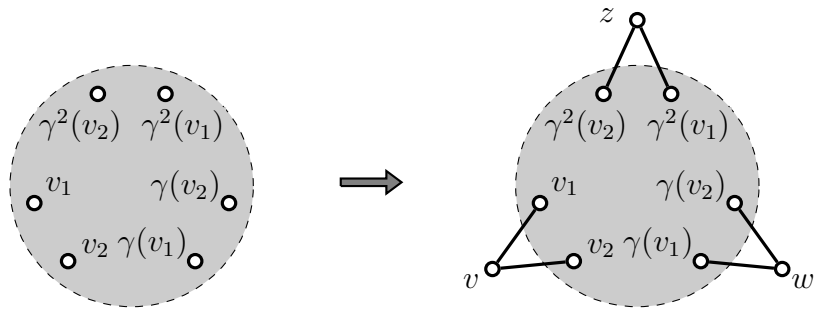

\begin{defin}
\label{C3vertex2extension}
\emph{Let $G$ be a graph, $\mathcal{C}_{3}=\{Id,C_{3},C_{3}^2\}$ be a symmetry group in dimension $2$, and $\Phi:\mathcal{C}_{3}\to \textrm{Aut}(G)$ be a homomorphism. Let $v_{1},v_{2}$ be two distinct vertices of $G$ and $v,w,z\notin V(G)$. Then the graph $\widehat{G}$ with $V(\widehat{G})=V(G)\cup \{v,w,z\}$ and $E(\widehat{G})=E(G)\cup \big\{\{v,v_{1}\},\{v,v_{2}\},\{w,\Phi(C_{3})(v_{1})\},\{w,\Phi(C_{3})(v_{2})\},\{z,\Phi(C_{3}^2)(v_{1})\},$ $\{z,\Phi(C_{3}^2)(v_{2})\}\big\}$ is called a $(\mathcal{C}_{3},\Phi)$ \emph{vertex addition (by $(v,w,z)$) of} $G$.}
\end{defin}

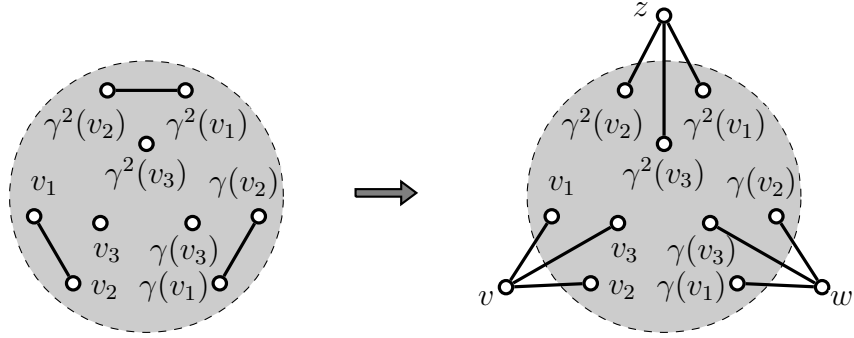
\begin{figure}[htp]
\begin{center}
\begin{tikzpicture}[very thick,scale=1]
\tikzstyle{every node}=[circle, draw=black, fill=white, inner sep=0pt, minimum width=5pt];
\filldraw[fill=black!20!white, draw=black, thin, dashed](0,0)circle(1.8cm);
\node (p1) at (70:1.5cm) {};
\node (p2) at (110:1.5cm) {};
\node (p3) at (190:1.5cm) {};
\node (p4) at (230:1.5cm) {};
\node (p5) at (310:1.5cm) {};
\node (p6) at (350:1.5cm) {};
\node (p7) at (90:0.7cm) {};
\node (p8) at (210:0.7cm) {};
\node (p9) at (330:0.7cm) {};
\node[rectangle,draw=black!20!white,fill=black!20!white](l1) at (50:1.25cm){$\gamma^2(v_1)$};
\node[rectangle,draw=black!20!white,fill=black!20!white](l2) at (130:1.25cm){$\gamma^2(v_2)$};
\node[rectangle,draw=black!20!white,fill=black!20!white](l3) at (174:1.35cm){$v_1$};
\node[rectangle,draw=black!20!white,fill=black!20!white](l4) at (246:1.35cm){$v_2$};
\node[rectangle,draw=black!20!white,fill=black!20!white](l5) at (287:1.25cm){$\gamma(v_1)$};
\node[rectangle,draw=black!20!white,fill=black!20!white](l6) at (368:1.3cm){$\gamma(v_2)$};
\node[rectangle,draw=black!20!white,fill=black!20!white](l7) at (90:0.3cm){$\gamma^2(v_3)$};
\node[rectangle,draw=black!20!white,fill=black!20!white](l8) at (235:0.9cm){$v_3$};
\node[rectangle,draw=black!20!white,fill=black!20!white](l9) at (305:0.88cm){$\gamma(v_3)$};
\draw(p1)--(p2);
\draw(p3)--(p4);
\draw(p5)--(p6);
\filldraw[fill=black!50!white, draw=black, thick]
    (2.75,0) -- (3.35,0) -- (3.35,-0.1) -- (3.55,0.05) -- (3.35,0.2) -- (3.35,0.1) -- (2.75,0.1) -- (2.75,0);
\end{tikzpicture}
\hspace{0.5cm}
\begin{tikzpicture}[very thick,scale=1]
\tikzstyle{every node}=[circle, draw=black, fill=white, inner sep=0pt, minimum width=5pt];
\filldraw[fill=black!20!white, draw=black, thin, dashed](0,0)circle(1.8cm);
\node (p1) at (70:1.5cm) {};
\node (p2) at (110:1.5cm) {};
\node (p3) at (190:1.5cm) {};
\node (p4) at (230:1.5cm) {};
\node (p5) at (310:1.5cm) {};
\node (p6) at (350:1.5cm) {};
\node (a7) at (90:2.4cm) {};
\node (a8) at (210:2.4cm) {};
\node (a9) at (330:2.4cm) {};
\node (p7) at (90:0.7cm) {};
\node (p8) at (210:0.7cm) {};
\node (p9) at (330:0.7cm) {};
\node[rectangle,draw=black!20!white,fill=black!20!white](l1) at (50:1.25cm){$\gamma^2(v_1)$};
\node[rectangle,draw=black!20!white,fill=black!20!white](l2) at (130:1.25cm){$\gamma^2(v_2)$};
\node[rectangle,draw=black!20!white,fill=black!20!white](l3) at (174:1.35cm){$v_1$};
\node[rectangle,draw=black!20!white,fill=black!20!white](l4) at (246:1.35cm){$v_2$};
\node[rectangle,draw=black!20!white,fill=black!20!white](l5) at (287:1.25cm){$\gamma(v_1)$};
\node[rectangle,draw=black!20!white,fill=black!20!white](l6) at (368:1.3cm){$\gamma(v_2)$};
\node[rectangle,draw=black!20!white,fill=black!20!white](l7) at (90:0.3cm){$\gamma^2(v_3)$};
\node[rectangle,draw=black!20!white,fill=black!20!white](l8) at (235:0.9cm){$v_3$};
\node[rectangle,draw=black!20!white,fill=black!20!white](l9) at (305:0.88cm){$\gamma(v_3)$};
\node[rectangle,draw=white,fill=white](l7) at (97:2.5cm){$z$};
\node[rectangle,draw=white,fill=white](l8) at (210:2.7cm){$v$};
\node[rectangle,draw=white,fill=white](l9) at (330:2.7cm){$w$};
\draw(a7)--(p1);
\draw(a7)--(p2);
\draw(a7)--(p7);
\draw(a8)--(p3);
\draw(a8)--(p4);
\draw(a8)--(p8);
\draw(a9)--(p5);
\draw(a9)--(p6);
\draw(a9)--(p9);
\end{tikzpicture}
\end{center}
\vspace{-0.3cm}
\caption{\emph{A $(\mathcal{C}_{3},\Phi)$ edge split of a graph $G$, where $\Phi(C_{3})=\gamma$ and $\Phi(C_{3}^2)=\gamma^2$.}}
\label{C3edspl}
\end{figure}

\begin{defin}
\label{C3edge2split}
\emph{Let $G$ be a graph, $\mathcal{C}_{3}=\{Id,C_{3},C_{3}^2\}$ be a symmetry group in dimension $2$, and $\Phi:\mathcal{C}_{3}\to \textrm{Aut}(G)$ be a homomorphism. Let $v_{1},v_{2},v_{3}$ be three distinct vertices of $G$ such that $\{v_{1},v_{2}\}\in E(G)$ and not both of $v_{1}$ and $v_{2}$ are fixed by $\Phi(C_{3})$ and let $v,w,z\notin V(G)$.
Then the graph $\widehat{G}$ with $V(\widehat{G})=V(G)\cup \{v,w,z\}$ and $E(\widehat{G})=\big(E(G)\setminus \big\{\{v_{1},v_{2}\},\{\Phi(C_{3})(v_{1}),\Phi(C_{3})(v_{2})\},\{\Phi(C_{3}^2)(v_{1}),\Phi(C_{3}^2)(v_{2})\}\big\}\big)\cup \big\{\{v,v_{i}\}|\, i=1,2,3\big\}\cup \big\{\{w,\Phi(C_{3})(v_{i})\}|\, i=1,2,3\big\}\cup \big\{\{z,\Phi(C_{3}^2)(v_{i})\}|\, i=1,2,3\big\}$ is called a $(\mathcal{C}_{3},\Phi)$ \emph{edge split (on $(\{v_{1},v_{2}\},\{\Phi(C_{3})(v_{1}),\Phi(C_{3})(v_{2})\},\{\Phi(C_{3}^2)(v_{1}),\Phi(C_{3}^2)(v_{2})\});(v,w,z)$) of $G$.}}
\end{defin}

\begin{figure}[htp]
\begin{center}
\begin{tikzpicture}[very thick,scale=1]
\tikzstyle{every node}=[circle, draw=black, fill=white, inner sep=0pt, minimum width=5pt];
\filldraw[fill=black!20!white, draw=black, thin, dashed](0,0)circle(1.3cm);
\node (p1) at (30:1cm) {};
\node (p2) at (150:1cm) {};
\node (p3) at (270:1cm) {};
\node[rectangle,draw=black!20!white,fill=black!20!white](l1) at (10:0.69cm){$\gamma(v_0)$};
\node[rectangle,draw=black!20!white,fill=black!20!white](l2) at (170:0.69cm){$\gamma^2(v_0)$};
\node[rectangle,draw=black!20!white,fill=black!20!white](l3) at (270:0.7cm){$v_0$};
\filldraw[fill=black!50!white, draw=black, thick]
    (2.75,0) -- (3.35,0) -- (3.35,-0.1) -- (3.55,0.05) -- (3.35,0.2) -- (3.35,0.1) -- (2.75,0.1) -- (2.75,0);
\end{tikzpicture}
\hspace{0.5cm}
\begin{tikzpicture}[very thick,scale=1]
\tikzstyle{every node}=[circle, draw=black, fill=white, inner sep=0pt, minimum width=5pt];
\filldraw[fill=black!20!white, draw=black, thin, dashed](0,0)circle(1.3cm);
\node (p1) at (30:1cm) {};
\node (p2) at (150:1cm) {};
\node (p3) at (270:1cm) {};
\node (p4) at (90:2.8cm) {};
\node (p5) at (210:2.8cm) {};
\node (p6) at (330:2.8cm) {};
\node[rectangle,draw=black!20!white,fill=black!20!white](l1) at (10:0.69cm){$\gamma(v_0)$};
\node[rectangle,draw=black!20!white,fill=black!20!white](l2) at (170:0.69cm){$\gamma^2(v_0)$};
\node[rectangle,draw=black!20!white,fill=black!20!white](l3) at (270:0.7cm){$v_0$};
\node[rectangle,draw=white,fill=white](l4) at (97:2.8cm){$z$};
\node[rectangle,draw=white,fill=white](l5) at (208:3.11cm){$v$};
\node[rectangle,draw=white,fill=white](l6) at (332:3.11cm){$w$};
\draw(p4)--(p5);
\draw(p5)--(p6);
\draw(p6)--(p4);
\draw(p4)--(p2);
\draw(p5)--(p3);
\draw(p6)--(p1);
\end{tikzpicture}
\end{center}
\vspace{-0.3cm}
\caption{\emph{A $(\mathcal{C}_{3},\Phi)$ $\Delta$ extension of a graph $G$, where $\Phi(C_{3})=\gamma$ and $\Phi(C_{3}^2)=\gamma^2$.}}
\label{C3delta}
\end{figure}
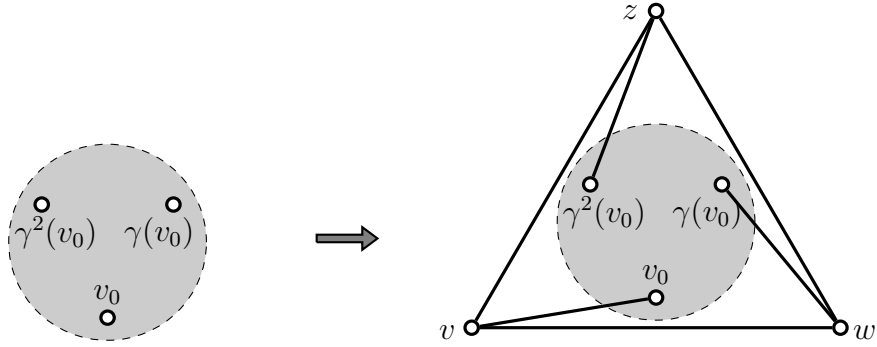

\begin{defin}
\label{C3delta2ext}
\emph{Let $G$ be a graph, $\mathcal{C}_{3}=\{Id,C_{3},C_{3}^2\}$ be a symmetry group in dimension $2$, and $\Phi:\mathcal{C}_{3}\to \textrm{Aut}(G)$ be a homomorphism. Let $v_{0}$ be a vertex of $G$ that is not fixed by $\Phi(C_{3})$ and let $v,w,z\notin V(G)$. Then the graph $\widehat{G}$ with $V(\widehat{G})=V(G)\cup \{v,w,z\}$ and $E(\widehat{G})=E(G)\cup \big\{\{v,w\},\{w,z\},\{z,v\},\{v,v_{0}\},\{w,\Phi(C_{3})(v_{0})\},\{z,\Phi(C_{3}^2)(v_{0})\}\big\}$ is called a $(\mathcal{C}_{3},\Phi)$ $\Delta$ \emph{extension (by $(v,w,z)$) of} $G$.}
\end{defin}

\begin{remark}
\label{neededforproper2}
\emph{Each of the constructions in Definitions \ref{C3vertex2extension}, \ref{C3edge2split}, and \ref{C3delta2ext} has the property that if the graph $G$ satisfies the Laman conditions, then so does $\widehat{G}$. This follows from Theorems \ref{lamantheorem} and \ref{henneberg} and the fact that we can obtain a $(\mathcal{C}_{3},\Phi)$ vertex addition of $G$ by a sequence of three vertex 2-additions, a $(\mathcal{C}_{3},\Phi)$ edge split of $G$ by a sequence of three edge 2-splits, and a $(\mathcal{C}_{3},\Phi)$ $\Delta$ extension of $G$ by a vertex 2-addition followed by two edge 2-splits.}
\end{remark}

In order to extend Crapo's Theorem to $\mathcal{C}_{3}$ we need the following symmetrized definition of a 3Tree2 partition.

\begin{defin}
\label{3T2C3}
\emph{Let $G$ be a graph, $\mathcal{C}_{3}=\{Id,C_{3},C_{3}^2\}$ be a symmetry group in dimension $2$, and $\Phi:\mathcal{C}_{3}\to \textrm{Aut}(G)$ be a homomorphism. A \emph{$(\mathcal{C}_{3},\Phi)$ 3Tree2 partition} of $G$ is a 3Tree2 partition $\{E(T_{0}),E(T_{1}),E(T_{2})\}$ of $G$ such that $\Phi(C_{3})(T_{i})= T_{i+1}$ for $i=0,1,2$, where the indices are added modulo 3.}
\end{defin}

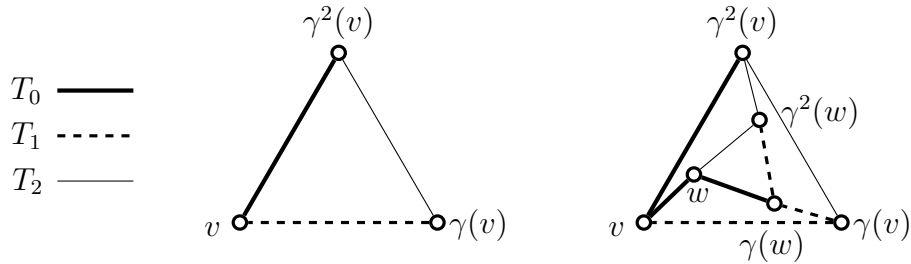
\begin{figure}[htp]
\begin{center}
\begin{tikzpicture}[very thick,scale=1]
\tikzstyle{every node}=[circle, draw=black, fill=white, inner sep=0pt, minimum width=5pt];
\draw [ultra thick] (-3.7,1) node[rectangle, draw=white,fill=white,left =5pt] {$T_0$} -- (-2.7,1);
\draw [dashed] (-3.7,0.4)  node[rectangle, draw=white,fill=white,left =5pt] {$T_1$} -- (-2.7,0.4);
\draw [thin] (-3.7,-0.2)  node[rectangle, draw=white,fill=white,left =5pt] {$T_2$} -- (-2.7,-0.2);
\node (p1) at (90:1.5cm) {};
\node (p2) at (210:1.5cm) {};
\node (p3) at (330:1.5cm) {};
\node[rectangle,draw=white,fill=white](l1) at (90:1.9cm){$\gamma^2(v)$};
\node[rectangle,draw=white,fill=white](l2) at (207:1.86cm){$v$};
\node[rectangle,draw=white,fill=white](l3) at (337:2cm){$\gamma(v)$};
\draw[ultra thick](p1)--(p2);
\draw[dashed](p3)--(p2);
\draw[thin](p1)--(p3);
\node [draw=white, fill=white] (a) at (270:1.6cm) {};
\end{tikzpicture}
\hspace{1cm}
\begin{tikzpicture}[very thick,scale=1]
\tikzstyle{every node}=[circle, draw=black, fill=white, inner sep=0pt, minimum width=5pt];
\node (b1) at (90:1.5cm) {};
\node (b2) at (210:1.5cm) {};
\node (b3) at (330:1.5cm) {};
\node(b4) at (70:0.65cm) {};
\node(b5) at (190:0.65cm) {};
\node(b6) at (310:0.65cm) {};
\node [rectangle,draw=white, fill=white] (w) at (213:0.7cm) {$w$};
\node [rectangle,draw=white, fill=white] (w2) at (290:1.13cm) {$\gamma(w)$};
\node [rectangle,draw=white, fill=white] (w3) at (35:1.25cm) {$\gamma^2(w)$};
\node [draw=white, fill=white] (b) at (270:1.6cm) {};
\draw [ultra thick](b1)--(b2);
\draw [thin](b1)--(b3);
\draw [dashed](b3)--(b2);
\draw [thin](b1)--(b4);
\draw [ultra thick](b2)--(b5);
\draw [dashed](b3)--(b6);
\draw [thin](b4)--(b5);
\draw [ultra thick](b5)--(b6);
\draw [dashed](b6)--(b4);
\node[rectangle,draw=white,fill=white](l1) at (90:1.9cm){$\gamma^2(v)$};
\node[rectangle,draw=white,fill=white](l2) at (207:1.86cm){$v$};
\node[rectangle,draw=white,fill=white](l3) at (337:2cm){$\gamma(v)$};
\end{tikzpicture}
\end{center}
\vspace{-0.3cm}
\caption{\emph{$(\mathcal{C}_{3},\Phi)$ 3Tree2 partitions of graphs, where $\Phi(C_{3})=\gamma$ and $\Phi(C_{3}^2)=\gamma^2$.}}
\label{C3te}
\end{figure}

\section{The main result}
\label{sec:main}

\begin{theorem}
\label{C3char}
Let $G$ be a graph with $|V(G)|\geq 3$, $\mathcal{C}_{3}=\{Id,C_{3},C_{3}^2\}$ be a symmetry group in dimension $2$, and $\Phi:\mathcal{C}_{3}\to \textrm{Aut}(G)$ be a homomorphism. The following are equivalent:
\begin{itemize}
\item[(i)] $\mathscr{R}_{(G,\mathcal{C}_{3},\Phi)}\neq \emptyset$ and $G$ is $(\mathcal{C}_{3},\Phi)$-generically isostatic;
\item[(ii)] $|E(G)|=2|V(G)|-3$, $|E(H)|\leq 2|V(H)|-3$ for all $H\subseteq G$ with $|V(H)|\geq 2$ (Laman conditions), and $j_{\Phi(C_{3})}=0$;
\item[(iii)] there exists a $(\mathcal{C}_{3},\Phi)$ construction sequence \begin{displaymath}(K_{3},\Phi_{0})=(G_{0},\Phi_{0}),(G_{1},\Phi_{1}),\ldots,(G_{k},\Phi_{k})=(G,\Phi)\end{displaymath} such that
    \begin{itemize}
    \item[(a)] $G_{i+1}$ is a $(\mathcal{C}_{3},\Phi_{i})$ vertex addition, a $(\mathcal{C}_{3},\Phi_{i})$ edge split, or a $(\mathcal{C}_{3},\Phi_{i})$ $\Delta$ extension of $G_{i}$ with $V(G_{i+1})=V(G_{i})\cup \{v_{i+1},w_{i+1},z_{i+1}\}$ for all $i=0,1,\ldots,k-1$;
    \item[(b)]
         $\Phi_{0}:\mathcal{C}_{3}\to \textrm{Aut}(K_{3})$ is a non-trivial homomorphism and for all $i=0,1,\ldots,k-1$, $\Phi_{i+1}:\mathcal{C}_{3}\to \textrm{Aut}(G_{i+1})$ is the homomorphism defined by $\Phi_{i+1}(x)|_{V(G_{i})}=\Phi_{i}(x)$ for all $x\in \mathcal{C}_{3}$ and $\Phi_{i+1}(C_{3})|_{\{v_{i+1},w_{i+1},z_{i+1}\}}=(v_{i+1}\,w_{i+1}\,z_{i+1})$;
    \end{itemize}
\item[(iv)] $G$ has a proper $(\mathcal{C}_{3},\Phi)$ 3Tree2 partition.
\end{itemize}
\end{theorem}

We break the proof of this result up into four Lemmas.

\begin{lemma}\label{3itoii} Let $G$ be a graph with $|V(G)|\geq 3$, $\mathcal{C}_{3}=\{Id,C_{3},C_{3}^2\}$ be a symmetry group in dimension $2$, and $\Phi:\mathcal{C}_{3}\to \textrm{Aut}(G)$ be a homomorphism. If $\mathscr{R}_{(G,\mathcal{C}_{3},\Phi)}\neq \emptyset$ and $G$ is $(\mathcal{C}_{3},\Phi)$-generically isostatic, then $G$ satisfies the Laman conditions and we have $j_{\Phi(C_{3})}=0$.
\end{lemma}
\textbf{Proof.} The result follows immediately from Laman's Theorem (Theorem \ref{lamantheorem}) and Theorem \ref{symmaxrule}. $\square$

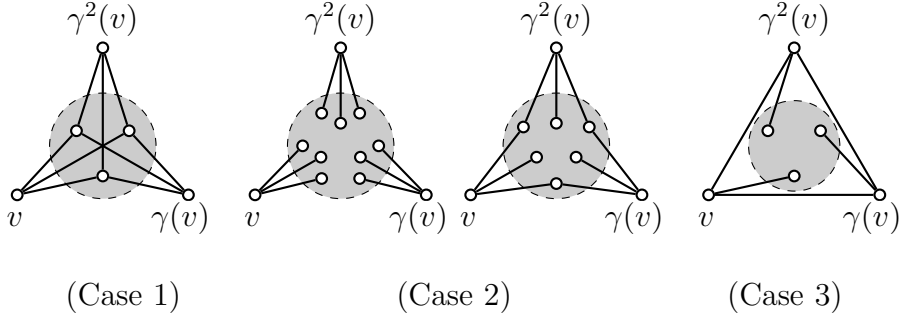
\begin{figure}[htp]
\begin{center}
\begin{tikzpicture}[thick,scale=1]
\tikzstyle{every node}=[circle, draw=black, fill=white, inner sep=0pt, minimum width=4pt];
\filldraw[fill=black!20!white, draw=black, thin, dashed](0,0)circle(0.7cm);
\node (p1) at (30:0.4cm) {};
\node (p2) at (150:0.4cm) {};
\node (p3) at (270:0.4cm) {};
\node (p4) at (90:1.3cm) {};
\node (p5) at (210:1.3cm) {};
\node (p6) at (330:1.3cm) {};
\draw(p1)--(p4);
\draw(p1)--(p5);
\draw(p1)--(p6);
\draw(p2)--(p4);
\draw(p2)--(p5);
\draw(p2)--(p6);
\draw(p3)--(p4);
\draw(p3)--(p5);
\draw(p3)--(p6);
\node[rectangle,draw=white,fill=white](l1) at (90:1.68cm){$\gamma^2(v)$};
\node[rectangle,draw=white,fill=white](l2) at (220:1.51cm){$v$};
\node[rectangle,draw=white,fill=white](l3) at (317:1.4cm){$\gamma(v)$};
\end{tikzpicture}
\hspace{0.15cm}
\begin{tikzpicture}[thick,scale=1]
\tikzstyle{every node}=[circle, draw=black, fill=white, inner sep=0pt, minimum width=4pt];
\filldraw[fill=black!20!white, draw=black, thin, dashed](0,0)circle(0.7cm);
\node (p1) at (90:0.3cm) {};
\node (p2) at (210:0.3cm) {};
\node (p3) at (330:0.3cm) {};
\node (p11) at (60:0.5cm) {};
\node (p12) at (120:0.5cm) {};
\node (p21) at (180:0.5cm) {};
\node (p22) at (240:0.5cm) {};
\node (p31) at (300:0.5cm) {};
\node (p32) at (360:0.5cm) {};
\node (p4) at (90:1.3cm) {};
\node (p5) at (210:1.3cm) {};
\node (p6) at (330:1.3cm) {};
\draw(p4)--(p1);
\draw(p4)--(p11);
\draw(p4)--(p12);
\draw(p5)--(p2);
\draw(p5)--(p21);
\draw(p5)--(p22);
\draw(p6)--(p3);
\draw(p6)--(p31);
\draw(p6)--(p32);
\node[rectangle,draw=white,fill=white](l1) at (90:1.68cm){$\gamma^2(v)$};
\node[rectangle,draw=white,fill=white](l2) at (220:1.51cm){$v$};
\node[rectangle,draw=white,fill=white](l3) at (317:1.4cm){$\gamma(v)$};
\end{tikzpicture}
\begin{tikzpicture}[thick,scale=1]
\tikzstyle{every node}=[circle, draw=black, fill=white, inner sep=0pt, minimum width=4pt];
\filldraw[fill=black!20!white, draw=black, thin, dashed](0,0)circle(0.7cm);
\node (p1) at (90:0.3cm) {};
\node (p2) at (210:0.3cm) {};
\node (p3) at (330:0.3cm) {};
\node (p11) at (30:0.5cm) {};
\node (p22) at (150:0.5cm) {};
\node (p33) at (270:0.5cm) {};
\node (p4) at (90:1.3cm) {};
\node (p5) at (210:1.3cm) {};
\node (p6) at (330:1.3cm) {};
\draw(p4)--(p1);
\draw(p4)--(p11);
\draw(p4)--(p22);
\draw(p5)--(p2);
\draw(p5)--(p33);
\draw(p5)--(p22);
\draw(p6)--(p3);
\draw(p6)--(p33);
\draw(p6)--(p11);
\node[rectangle,draw=white,fill=white](l1) at (90:1.68cm){$\gamma^2(v)$};
\node[rectangle,draw=white,fill=white](l2) at (220:1.51cm){$v$};
\node[rectangle,draw=white,fill=white](l3) at (317:1.4cm){$\gamma(v)$};
\end{tikzpicture}
\hspace{0.15cm}
\begin{tikzpicture}[thick,scale=1]
\tikzstyle{every node}=[circle, draw=black, fill=white, inner sep=0pt, minimum width=4pt];
\filldraw[fill=black!20!white, draw=black, thin, dashed](0,0)circle(0.6cm);
\node (p1) at (30:0.4cm) {};
\node (p2) at (150:0.4cm) {};
\node (p3) at (270:0.4cm) {};
\node (p4) at (90:1.3cm) {};
\node (p5) at (210:1.3cm) {};
\node (p6) at (330:1.3cm) {};
\node[rectangle,draw=white,fill=white](l1) at (90:1.68cm){$\gamma^2(v)$};
\node[rectangle,draw=white,fill=white](l2) at (220:1.51cm){$v$};
\node[rectangle,draw=white,fill=white](l3) at (317:1.4cm){$\gamma(v)$};
\draw(p4)--(p5);
\draw(p6)--(p5);
\draw(p4)--(p6);
\draw(p4)--(p2);
\draw(p5)--(p3);
\draw(p6)--(p1);
\end{tikzpicture}
\begin{tikzpicture}[thick,scale=1]
\tikzstyle{every node}=[circle, draw=black, fill=white, inner sep=0pt, minimum width=4pt];
\node [draw=white, fill=white] (a) at (1.15,0) {(Case 1)};
\node [draw=white, fill=white] (a) at (5.5,0) {(Case 2)};
\node [draw=white, fill=white] (a) at (9.85,0) {(Case 3)};
\end{tikzpicture}
\end{center}
\vspace{-0.3cm}
\caption{\emph{If a graph $G$ satisfies the conditions in Theorem \ref{C3char} $(ii)$ and has a vertex $v$ of valence $3$, then $G$ is a graph of one of the types depicted above.}}
\label{3casesC3}
\end{figure}

\begin{lemma} \label{3iitoiii} Let $G$ be a graph with $|V(G)|\geq 3$, $\mathcal{C}_{3}=\{Id,C_{3},C_{3}^2\}$ be a symmetry group in dimension $2$, and $\Phi:\mathcal{C}_{3}\to \textrm{Aut}(G)$ be a homomorphism. If $G$ satisfies the Laman conditions and we also have $j_{\Phi(C_{3})}=0$, then there exists a $(\mathcal{C}_{3},\Phi)$ construction sequence for $G$.
\end{lemma}
\textbf{Proof.} We employ induction on $|V(G)|$. Note first that if for a graph $G$, there exists a homomorphism $\Phi:\mathcal{C}_{3}\to \textrm{Aut}(G)$ such that $j_{\Phi(C_{3})}=0$, then $|V(G)|\equiv 0 \pmod{3}$. The only graph with three vertices that satisfies the Laman conditions is the graph $K_{3}$ and if $\Phi:\mathcal{C}_{3}\to \textrm{Aut}(K_{3})$ is a homomorphism such that $j_{\Phi(C_{3})}=0$, then $\Phi$ is clearly a non-trivial homomorphism. This proves the base case.\\\indent So we let $n\geq 3$ and we assume that the result holds for all graphs with $n$ or fewer than $n$ vertices.\\\indent Let $G$ be a graph with $|V(G)|=n+3$ that satisfies the Laman conditions and suppose $j_{\Phi(C_{3})}=0$ for a homomorphism $\Phi:\mathcal{C}_{3}\to \textrm{Aut}(G)$. In the following, we denote $\Phi(C_{3})$ by $\gamma$ and $\Phi(C_{3}^2)$ by $\gamma^2$.\\\indent Since $G$ satisfies the Laman conditions, it is easy to verify that $G$ has a vertex of valence $2$ or $3$ (see \cite{graver, gss, BS4}, for example).\\\indent
We assume first that $G$ has a vertex $v$ of valence $2$, say $N_{G}(v)=\{v_{1},v_{2}\}$. Note that $v, \gamma(v)$ and $\gamma^2(v)$ are three distinct vertices of $G$, because $j_{\gamma}=0$. Suppose two of these vertices are adjacent, wlog $\{v,\gamma(v)\}\in E(G)$. Then $\{\gamma(v),\gamma^2(v)\},\{\gamma^2(v),v\}\in E(G)$, because $\gamma\in \textrm{Aut}(G)$. Let $G'=G-\{v, \gamma(v),\gamma^2(v)\}$. Then \begin{displaymath}|E(G')|=|E(G)|-3=2|V(G)|-6=2|V(G')|\textrm{.}\end{displaymath} Since $|V(G)|\geq 6$, we have $|V(G')|\geq 3$, and hence $G'$ violates the Laman conditions, a contradiction.\\\indent Therefore, $\{v, \gamma(v),\gamma^2(v)\}$ is an independent subset of $V(G)$, which says that the six edges $\{v,v_{i}\},\{\gamma(v),\gamma(v_{i})\},\{\gamma^2(v),\gamma^2(v_{i})\}$, $i=1,2$, are all pairwise distinct. Thus, \begin{displaymath}|E(G')|=|E(G)|-6=2|V(G)|-9=2|V(G')|-3\textrm{.}\end{displaymath} Also, for $H\subseteq G'$ with $|V(H)|\geq 2$, we have $H\subseteq G$, and hence \begin{displaymath}|E(H)|\leq 2|V(H)|-3\textrm{.}\end{displaymath} Therefore, $G'$ satisfies the Laman conditions.\\\indent Let $\Phi':\mathcal{C}_{3}\to \textrm{Aut}(G')$ be the homomorphism with $\Phi'(x)=\Phi(x)|_{V(G')}$ for all $x\in \mathcal{C}_{3}$. Then we have $j_{\Phi'(C_{3})}=0$, and hence, by the induction hypothesis, there exists a sequence \begin{displaymath}(K_{3},\Phi_{0})=(G_{0},\Phi_{0}),(G_{1},\Phi_{1}),\ldots,(G_{k},\Phi_{k})=(G',\Phi')\end{displaymath} satisfying the conditions in Theorem \ref{C3char} $(iii)$. Since $G$ is a $(\mathcal{C}_{3},\Phi')$ vertex addition of $G'$ with $V(G)=V(G')\cup\{v, \gamma(v),\gamma^2(v)\}$, \begin{displaymath}(K_{3},\Phi_{0})=(G_{0},\Phi_{0}),(G_{1},\Phi_{1}),\ldots,(G',\Phi'),(G,\Phi)\end{displaymath} is a sequence with the desired properties. \\\indent Suppose now that $G$ has a vertex $v$ of valence $3$, say $N_{G}(v)=\{v_{1},v_{2},v_{3}\}$, and no vertex of valence $2$. Note that $v, \gamma(v)$ and $\gamma^2(v)$ are again three distinct vertices of $G$, as are $v_{i}, \gamma(v_{i})$ and $\gamma^2(v_{i})$ for each $i=1,2,3$, because $j_{\gamma}=0$. We need to consider the following three cases (see also Figure \ref{3casesC3}):

\begin{itemize}
\item[]
\begin{itemize}
\item[ \textbf{Case 1:}] $\{v, \gamma(v),\gamma^2(v)\}$ is an independent subset of $V(G)$ and all three of these vertices share a common neighbor, say wlog $v_{1}$. Since $\gamma\in \textrm{Aut}(G)$, this says that each of $v, \gamma(v)$ and $\gamma^2(v)$ has the same neighbors, namely $v_{1}, \gamma(v_{1})$ and $\gamma^2(v_{1})$.
\item[\textbf{Case 2:}] $\{v, \gamma(v),\gamma^2(v)\}$ is an independent subset of $V(G)$ and there is no vertex in $G$ that is adjacent to all three of these vertices.
\item[\textbf{Case 3:}] $\{v, \gamma(v),\gamma^2(v)\}$ is not independent in $G$ and hence $\{v,\gamma(v)\},\{\gamma(v),\gamma^2(v)\},\{\gamma^2(v),v\}\in E(G)$.
\end{itemize}
\end{itemize}

\textbf{Case 1:} By Theorems \ref{lamantheorem} and \ref{edgespl}, there exists a pair $\{a,b\}$ of vertices in $\{v_{1},\gamma(v_{1}),\gamma^2(v_{1})\}$
such that $G-\{v\}+\big\{\{a,b\}\big\}$ satisfies the Laman conditions. By the same argument, applied two more times, it follows that the graph $\widetilde{G}= G-\big\{\{v, \gamma(v), \gamma^2(v)\}\big\}+ \big\{\{v_{1},\gamma(v_{1})\},\{\gamma(v_{1}),\gamma^2(v_{1})\},\{\gamma^2(v_{1}),v_{1}\}\big\}$ also satisfies the Laman conditions.
\\\indent Further, if we define $\widetilde{\Phi}$ by $\widetilde{\Phi}(x)=\Phi(x)|_{V(\widetilde{G})}$ for all $x\in \mathcal{C}_{3}$, then $\widetilde{\Phi}(x)\in \textrm{Aut}(\widetilde{G})$ for all $x\in \mathcal{C}_{3}$ and $\widetilde{\Phi}:\mathcal{C}_{3}\to \textrm{Aut}(\widetilde{G})$ is a homomorphism. Since we clearly also have $j_{\widetilde{\Phi}(C_{3})}=0$, it follows from the induction hypothesis that there exists a sequence \begin{displaymath}(K_{3},\Phi_{0})=(G_{0},\Phi_{0}),(G_{1},\Phi_{1}),\ldots,(G_{k},\Phi_{k})=(\widetilde{G},\widetilde{\Phi})
\end{displaymath} satisfying the conditions in Theorem \ref{C3char} $(iii)$. Since $G$ is a $(\mathcal{C}_{3},\widetilde{\Phi})$ edge split of $\widetilde{G}$ with $V(G)=V(\widetilde{G})\cup\{v, \gamma(v),\gamma^2(v)\}$, \begin{displaymath}(K_{3},\Phi_{0})=(G_{0},\Phi_{0}),(G_{1},\Phi_{1}),\ldots,(\widetilde{G},\widetilde{\Phi}),(G,\Phi)\end{displaymath} is a sequence with the desired properties.

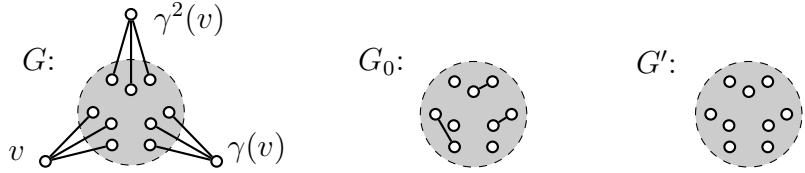
\begin{figure}[htp]
\begin{center}
\begin{tikzpicture}[thick,scale=1]
\tikzstyle{every node}=[circle, draw=black, fill=white, inner sep=0pt, minimum width=4pt];
\filldraw[fill=black!20!white, draw=black, thin, dashed](0,0)circle(0.7cm);
\node (p1) at (90:0.3cm) {};
\node (p2) at (210:0.3cm) {};
\node (p3) at (330:0.3cm) {};
\node (p11) at (60:0.5cm) {};
\node (p12) at (120:0.5cm) {};
\node (p21) at (180:0.5cm) {};
\node (p22) at (240:0.5cm) {};
\node (p31) at (300:0.5cm) {};
\node (p32) at (360:0.5cm) {};
\node (p4) at (90:1.3cm) {};
\node (p5) at (210:1.3cm) {};
\node (p6) at (330:1.3cm) {};
\draw(p4)--(p1);
\draw(p4)--(p11);
\draw(p4)--(p12);
\draw(p5)--(p2);
\draw(p5)--(p21);
\draw(p5)--(p22);
\draw(p6)--(p3);
\draw(p6)--(p31);
\draw(p6)--(p32);
\node[rectangle,draw=white,fill=white](l1) at (58:1.45cm){$\gamma^2(v)$};
\node[rectangle,draw=white,fill=white](l2) at (200:1.6cm){$v$};
\node[rectangle,draw=white,fill=white](l3) at (344:1.72cm){$\gamma(v)$};
\node [draw=white, fill=white] (a) at (150:1.4cm) {$G$:};
\end{tikzpicture}
\hspace{0.6cm}
\begin{tikzpicture}[thick,scale=1]
\tikzstyle{every node}=[circle, draw=black, fill=white, inner sep=0pt, minimum width=4pt];
\filldraw[fill=black!20!white, draw=black, thin, dashed](0,0)circle(0.7cm);
\node (p1) at (90:0.3cm) {};
\node (p2) at (210:0.3cm) {};
\node (p3) at (330:0.3cm) {};
\node (p11) at (60:0.5cm) {};
\node (p12) at (120:0.5cm) {};
\node (p21) at (180:0.5cm) {};
\node (p22) at (240:0.5cm) {};
\node (p31) at (300:0.5cm) {};
\node (p32) at (360:0.5cm) {};
\draw(p11)--(p1);
\draw(p21)--(p22);
\draw(p3)--(p32);
\node [draw=white, fill=white] (b) at (150:1.4cm) {$G_0$:};
\end{tikzpicture}
\hspace{1.1cm}
\begin{tikzpicture}[thick,scale=1]
\tikzstyle{every node}=[circle, draw=black, fill=white, inner sep=0pt, minimum width=4pt];
\filldraw[fill=black!20!white, draw=black, thin, dashed](0,0)circle(0.7cm);
\node (p1) at (90:0.3cm) {};
\node (p2) at (210:0.3cm) {};
\node (p3) at (330:0.3cm) {};
\node (p11) at (60:0.5cm) {};
\node (p12) at (120:0.5cm) {};
\node (p21) at (180:0.5cm) {};
\node (p22) at (240:0.5cm) {};
\node (p31) at (300:0.5cm) {};
\node (p32) at (360:0.5cm) {};
\node [draw=white, fill=white] (c) at (150:1.4cm) {$G'$:};
\end{tikzpicture}
\end{center}
\vspace{-0.3cm}
\caption{\emph{The graphs $G$, $G_0$ and $G'$ in Case 2 of the proof of Lemma \ref{3iitoiii}.}}
\label{figG}
\end{figure}

\textbf{Case 2:} By Theorems \ref{lamantheorem} and \ref{edgespl}, there exists $\{i_2,j_2\}\subseteq \{1,2,3\}$ such that $G_2=G-\{\gamma^2(v)\}+\big\{\{\gamma^2(v_{i_2}),\gamma^2(v_{j_2})\}\big\}$ satisfies the Laman conditions. By the same argument, there exist $\{i_1,j_1\}\subseteq \{1,2,3\}$ and $\{i_0,j_0\}\subseteq \{1,2,3\}$ such that both $G_1=G_2-\{\gamma(v)\}+\big\{\{\gamma(v_{i_1}),\gamma(v_{j_1})\}\big\}$ and $G_0=G_1-\{v\}+\big\{\{v_{i_0},v_{j_0}\}\big\}$ also satisfy the Laman conditions. We assume wlog that $\{i_0,j_0\}=\{1,2\}$. Then for every subgraph $H$ of $G_0-\big\{\{v_{1},v_{2}\}\big\}$ with $v_{1},v_{2}\in V(H)$ we have $|E(H)|\leq 2|V(H)|-4$. Moreover, for every subgraph $H$ of $G'=G-\{v, \gamma(v),\gamma^2(v)\}$ with $v_{1},v_{2}\in V(H)$ we also have $|E(H)|\leq 2|V(H)|-4$, because $G'$ is obtained from $G_0-\big\{\{v_{1},v_{2}\}\big\}$ by deleting the edges $\{\gamma(v_{i_1}),\gamma(v_{j_1})\}$ and $\{\gamma^2(v_{i_2}),\gamma^2(v_{j_2})\}$ (see also Figure \ref{figG}).

Since $G'$ is invariant under $\gamma$, every subgraph $H$ of $G'$ with $\gamma(v_{1}),\gamma(v_{2})\in V(H)$ or $\gamma^2(v_{1}),\gamma^2(v_{2})\in V(H)$ also satisfies $|E(H)|\leq 2|V(H)|-4$. Note that $\{v_{1},v_{2}\},\{\gamma(v_{1}),\gamma(v_{2})\}$ and $\{\gamma^2(v_{1}),\gamma^2(v_{2})\}$ are three distinct pairs of vertices (though not edges, by the above counts) of $G$, as the following argument shows.\\\indent Suppose $\{v_{1},v_{2}\}=\{\gamma(v_{1}),\gamma(v_{2})\}$. Then $v_{1}=\gamma(v_{2})$ and $v_{2}=\gamma(v_{1})$, because $G$ satisfies $j_{\gamma}=0$. Therefore, $\gamma(v_{1})=\gamma^2(v_{2})$, and hence $v_{2}=\gamma^2(v_{2})$, contradicting $j_{\gamma}=0$. Similarly, $\{\gamma(v_{1}),\gamma(v_{2})\}\neq \{\gamma^2(v_{1}),\gamma^2(v_{2})\}$ and $\{v_{1},v_{2}\}\neq \{\gamma^2(v_{1}),\gamma^2(v_{2})\}$.\\\indent We claim that $\widetilde{G}= G'+ \big\{\{v_{1},v_{2}\},\{\gamma(v_{1}),\gamma(v_{2})\},\{\gamma^2(v_{1}),\gamma^2(v_{2})\}\big\}$ satisfies the Laman conditions. We clearly have \begin{displaymath}|E(\widetilde{G})|=|E(G')|+3=|E(G)|-6=2|V(G)|-9=2|V(\widetilde{G})|-3\textrm{.}\end{displaymath} \indent Suppose there exists a subgraph $H$ of $G'$ with $v_{1},v_{2},\gamma(v_{1}),\gamma(v_{2})\in V(H)$ and $|E(H)|= 2|V(H)|-4$. Then there also exists $\gamma(H)\subseteq G'$ with $\gamma(v_{1}),\gamma(v_{2}),\gamma^2(v_{1}),\gamma^2(v_{2})\in V\big(\gamma(H)\big)$ and $|E\big(\gamma(H)\big)|= 2|V\big(\gamma(H)\big)|-4$, as well as $\gamma^2(H)\subseteq G'$ with $\gamma^2(v_{1}),\gamma^2(v_{2}),v_{1},v_{2}\in V\big(\gamma^2(H)\big)$ and $|E\big(\gamma^2(H)\big)|= 2|V\big(\gamma^2(H)\big)|-4$, because $G'$ is invariant under $\gamma$. Let $H'=H\cup \gamma(H)$. Then
\begin{eqnarray}
|E(H')|& = & |E(H)|+|E\big(\gamma(H)\big)|-|E\big(H\cap \gamma(H)\big)|\nonumber\\
 & \geq & 2|V(H)|-4+ 2|V\big(\gamma(H)\big)|-4-(2|V\big(H\cap \gamma(H)\big)|-4)\nonumber\\
 & = & 2|V(H')|-4 \textrm{,}\nonumber
\end{eqnarray}
because  $H\cap \gamma(H)$ is a subgraph of $G'$ with $\gamma(v_{1}),\gamma(v_{2})\in V\big(H\cap \gamma(H)\big)$. Since $H'$ is also a subgraph of $G'$ with $\gamma(v_{1}),\gamma(v_{2})\in V(H')$, it follows that \begin{displaymath}|E(H')|=2|V(H')|-4\textrm{.} \end{displaymath} Similarly, it can be shown that $H''=H'\cup \gamma^2(H)$ satisfies \begin{displaymath}|E(H'')|=2|V(H'')|-4\textrm{,} \end{displaymath} because $H'\cap \gamma^2(H)$ is a subgraph of $G'$ with $v_{1},v_{2}\in V\big(H'\cap \gamma^2(H)\big)$. However, $H''$ is invariant under $\gamma$ and satisfies $j_{\gamma|_{V(H'')}}=0$, so that $|V(H'')|\equiv 0 \pmod{3}$ and $|E(H'')|\equiv 0 \pmod{3}$, contradicting the count $|E(H'')|=2|V(H'')|-4$.\\\indent Therefore, every subgraph $H$ of $G'$ with $v_{1},v_{2},\gamma(v_{1}),\gamma(v_{2})\in V(H)$, $\gamma(v_{1}),\gamma(v_{2}),\gamma^2(v_{1}),\gamma^2(v_{2})\in V(H)$, or $\gamma^2(v_{1}),\gamma^2(v_{2}),v_{1},v_{2}\in V(H)$ satisfies $|E(H)|\leq 2|V(H)|-5$.\\\indent
It is now only left to show that for every subgraph $H$ of $G'$ with $v_{1},v_{2},\gamma(v_{1}),\gamma(v_{2}),\gamma^2(v_{1}),\gamma^2(v_{2})\in V(H)$, we have $|E(H)|\leq 2|V(H)|-6$. Suppose to the contrary that there exists a subgraph $H$ of $G'$ with $v_{1},v_{2},\gamma(v_{1}),\gamma(v_{2}),\gamma^2(v_{1}),\gamma^2(v_{2})\in V(H)$ and $|E(H)|= 2|V(H)|-5$. Then there also exist $\gamma(H)\subseteq G'$ and $\gamma^2(H)\subseteq G'$ with the same properties, because $G'$ is invariant under $\gamma$. Let $H'=H\cup \gamma(H)$. Then
\begin{eqnarray}
|E(H')|& = & |E(H)|+|E\big(\gamma(H)\big)|-|E\big(H\cap \gamma(H)\big)|\nonumber\\
 & \geq & 2|V(H)|-5+ 2|V\big(\gamma(H)\big)|-5-(2|V\big(H\cap \gamma(H)\big)|-5)\nonumber\\
 & = & 2|V(H')|-5 \textrm{,}\nonumber
\end{eqnarray}
because $H\cap \gamma(H)$ is a subgraph of $G'$ with $v_{1},v_{2},\gamma(v_{1}),\gamma(v_{2})\in V\big(H\cap \gamma(H)\big)$. Since $H'$ is also a subgraph of $G'$ with $v_{1},v_{2},\gamma(v_{1}),\gamma(v_{2})\in V(H')$, it follows that \begin{displaymath}|E(H')|= 2|V(H')|-5 \textrm{.} \end{displaymath} Similarly, it can be shown that $H''=H'\cup \gamma^2(H)$ satisfies \begin{displaymath}|E(H'')|=2|V(H'')|-5\textrm{,} \end{displaymath} because $H'\cap \gamma^2(H)$ is a subgraph of $G'$ with $v_{1},v_{2},\gamma(v_{1}),\gamma(v_{2})\in V\big(H'\cap \gamma^2(H)\big)$. However, $H''$ is invariant under $\gamma$ and we have $j_{\gamma|_{V(H'')}}=0$, so that $|V(H'')|\equiv 0 \pmod{3}$ and $|E(H'')|\equiv 0 \pmod{3}$, contradicting the count $|E(H'')|=2|V(H'')|-5$.\\\indent Thus, $\widetilde{G}= G'+ \big\{\{v_{1},v_{2}\},\{\gamma(v_{1}),\gamma(v_{2})\},\{\gamma^2(v_{1}),\gamma^2(v_{2})\}\big\}$ indeed satisfies the Laman conditions.\\\indent Further, if we define $\widetilde{\Phi}$ by $\widetilde{\Phi}(x)=\Phi(x)|_{V(\widetilde{G})}$ for all $x\in \mathcal{C}_{3}$, then $\widetilde{\Phi}(x)\in \textrm{Aut}(\widetilde{G})$ for all $x\in \mathcal{C}_{3}$ and $\widetilde{\Phi}:\mathcal{C}_{3}\to \textrm{Aut}(\widetilde{G})$ is a homomorphism. Since we also have $j_{\widetilde{\Phi}(C_{3})}=0$, it follows from the induction hypothesis that there exists a sequence \begin{displaymath}(K_{3},\Phi_{0})=(G_{0},\Phi_{0}),(G_{1},\Phi_{1}),\ldots,(G_{k},\Phi_{k})=(\widetilde{G},\widetilde{\Phi})
\end{displaymath} satisfying the conditions in Theorem \ref{C3char} $(iii)$. Since $G$ is a $(\mathcal{C}_{3},\widetilde{\Phi})$ edge split of $\widetilde{G}$ with $V(G)=V(\widetilde{G})\cup\{v, \gamma(v),\gamma^2(v)\}$, \begin{displaymath}(K_{3},\Phi_{0})=(G_{0},\Phi_{0}),(G_{1},\Phi_{1}),\ldots,(\widetilde{G},\widetilde{\Phi}),(G,\Phi)\end{displaymath} is a sequence with the desired properties.

\textbf{Case 3:} Note that $G'=G-\{v, \gamma(v),\gamma^2(v)\}$ satisfies \begin{displaymath}|E(G')|=|E(G)|-6=2|V(G)|-9=2|V(G')|-3\textrm{.}\end{displaymath} Also, for $H\subseteq G'$ with $|V(H)|\geq 2$, we have $H\subseteq G$, and hence \begin{displaymath}|E(H)|\leq 2|V(H)|-3\textrm{,}\end{displaymath} so that $G'$ satisfies the Laman conditions.\\\indent If we define $\Phi'$ by $\Phi'(x)=\Phi(x)|_{V(G')}$ for all $x\in \mathcal{C}_{3}$, then $\Phi'(x)\in \textrm{Aut}(G')$ for all $x\in \mathcal{C}_{3}$ and $\Phi':\mathcal{C}_{3}\to \textrm{Aut}(G')$ is a homomorphism. Since we also have $j_{\Phi'(C_{3})}=0$, it follows from the induction hypothesis that there exists a sequence \begin{displaymath}(K_{3},\Phi_{0})=(G_{0},\Phi_{0}),(G_{1},\Phi_{1}),\ldots,(G_{k},\Phi_{k})=(G',\Phi')\end{displaymath} satisfying the conditions in Theorem \ref{C3char} $(iii)$. Since $G$ is a $(\mathcal{C}_{3},\Phi')$ $\Delta$ extension of $G'$ with $V(G)=V(G')\cup\{v, \gamma(v),\gamma^2(v)\}$, \begin{displaymath}(K_{3},\Phi_{0})=(G_{0},\Phi_{0}),(G_{1},\Phi_{1}),\ldots,(G',\Phi'),(G,\Phi)\end{displaymath} is a sequence with the desired properties. $\square$

\begin{lemma}\label{3iiitoiv} Let $G$ be a graph with $|V(G)|\geq 3$, $\mathcal{C}_{3}=\{Id,C_{3},C_{3}^2\}$ be a symmetry group in dimension $2$, and $\Phi:\mathcal{C}_{3}\to \textrm{Aut}(G)$ be a homomorphism. If there exists a $(\mathcal{C}_{3},\Phi)$ construction sequence for $G$, then $G$ has a proper $(\mathcal{C}_{3},\Phi)$ 3Tree2 partition.
\end{lemma}
\textbf{Proof.} We proceed by induction on $|V(G)|$. Let $V(K_{3})=\{v_{1},v_{2},v_{3}\}$ and wlog let $\Phi:\mathcal{C}_{3}\to \textrm{Aut}(K_{3})$ be the homomorphism defined by $\Phi(C_{3})=(v_{1}\,v_{2}\,v_{3})$. Then $K_{3}$ has the proper $(\mathcal{C}_{3},\Phi)$ 3Tree2 partition $\{E(T_{0}),E(T_{1}),E(T_{2})\}$, where $T_{0}=\langle\{v_{1},v_{2}\}\rangle$, $T_{1}=\langle\{v_{2},v_{3}\}\rangle$ and $T_{2}=\langle\{v_{3},v_{1}\}\rangle$.\\\indent Assume, then, that the result holds for all graphs with $n$ or fewer than $n$ vertices, where $n\geq 3$.\\\indent
Let $G$ be a graph with $|V(G)|=n+3$ and let $\Phi:\mathcal{C}_{3}\to \textrm{Aut}(G)$ be a homomorphism such that there exists a $(\mathcal{C}_{3},\Phi)$ construction sequence \begin{displaymath}(K_{3},\Phi_{0})=(G_{0},\Phi_{0}),(G_{1},\Phi_{1}),\ldots,(G_{k},\Phi_{k})=(G,\Phi)\end{displaymath} satisfying the conditions in Theorem \ref{C3char} $(iii)$. By Remark \ref{neededforproper2}, $G$ satisfies the Laman conditions, and hence, by Remark \ref{3t2remk}, any 3Tree2 partition of $G$ must be proper. Therefore,  it suffices to show that $G$ has some $(\mathcal{C}_{3},\Phi)$ 3Tree2 partition. We let $\Phi(C_{3})=\gamma$ and $\Phi(C_{3}^2)=\gamma^2$.\\\indent By the induction hypothesis, $G_{k-1}$ has a $(\mathcal{C}_{3},\Phi_{k-1})$ 3Tree2 partition $\big\{E\big(T^{(k-1)}_{0}\big),E\big(T^{(k-1)}_{1}\big),E\big(T^{(k-1)}_{2}\big)\big\}$. In the following, we compute the indices $i$ of the trees $T^{(k-1)}_{i}$ modulo 3.\\\indent Suppose first that $G$ is a $(\mathcal{C}_{3},\Phi_{k-1})$ vertex addition by $(v\,w\,z)$ of $G_{k-1}$ with $N_{G}(v)=\{v_{1},v_{2}\}$, where $w=\gamma(v)$ and $z=\gamma^2(v)$. Since $\Phi_{k-1}(C_{3})=\gamma|_{V(G_{k-1})}$ we have $N_{G}(w)=\{\gamma(v_{1}),\gamma(v_{2})\}$ and $N_{G}(z)=\{\gamma^2(v_{1}),\gamma^2(v_{2})\}$. Note that both $v_{1}$ and $v_{2}$ belong to exactly two of the trees $T^{(k-1)}_{i}$. Therefore, there exists $l\in\{0,1,2\}$ such that $v_{1}\in V\big(T^{(k-1)}_{l}\big)$ and $v_{2}\in V\big(T^{(k-1)}_{l+1}\big)$. It follows that $\gamma(v_{1})\in V\big(T^{(k-1)}_{l+1}\big)$, $\gamma^2(v_{1}),\gamma(v_{2})\in V\big(T^{(k-1)}_{l+2}\big)$ and $\gamma^2(v_{2})\in V\big(T^{(k-1)}_{l}\big)$. So, if we define $T^{(k)}_{l}$ to be the tree with \begin{eqnarray}V\big(T^{(k)}_{l}\big)&=&V\big(T^{(k-1)}_{l}\big)\cup\{v,z\}\nonumber\\ E\big(T^{(k)}_{l}\big)&=&E\big(T^{(k-1)}_{l}\big)\cup\big\{\{v,v_{1}\},\{z,\gamma^2(v_{2})\}\big\}\textrm{,}\nonumber\end{eqnarray} $T^{(k)}_{l+1}$ to be the tree with
\begin{eqnarray}V\big(T^{(k)}_{l+1}\big)&=&V\big(T^{(k-1)}_{l+1}\big)\cup\{v,w\}\nonumber\\ E\big(T^{(k)}_{l+1}\big)&=&E\big(T^{(k-1)}_{l+1}\big)\cup\big\{\{v,v_{2}\},\{w,\gamma(v_{1})\}\big\}\textrm{,}\nonumber\end{eqnarray}
and $T^{(k)}_{l+2}$ to be the tree with
\begin{eqnarray}V\big(T^{(k)}_{l+2}\big)&=&V\big(T^{(k-1)}_{l+2}\big)\cup\{w,z\}\nonumber\\  E\big(T^{(k)}_{l+2}\big)&=&E\big(T^{(k-1)}_{l+2}\big)\cup\big\{\{w,\gamma(v_{2})\},\{z,\gamma^2(v_{1})\}\big\}\textrm{,}\nonumber
\end{eqnarray} then $\big\{E\big(T^{(k)}_{0}\big),E\big(T^{(k)}_{1}\big),E\big(T^{(k)}_{2}\big)\big\}$ is a $(\mathcal{C}_{3},\Phi)$ 3Tree2 partition of $G$.

\begin{figure}[htp]
\begin{center}
\begin{tikzpicture}[very thick,scale=1]
\tikzstyle{every node}=[circle, draw=black, fill=white, inner sep=0pt, minimum width=5pt];
\filldraw[fill=black!20!white, draw=black, thin, dashed](0,0)circle(1.6cm);
\node[rectangle,draw=black!20!white,fill=black!20!white](l5) at (290:0.81cm){$\gamma(v_1)$};
\node[rectangle,draw=black!20!white,fill=black!20!white](l6) at (360:0.71cm){$\gamma(v_2)$};
\node (p1) at (70:1.3cm) {};
\node (p2) at (110:1.3cm) {};
\node (p3) at (190:1.3cm) {};
\node (p4) at (230:1.3cm) {};
\node (p5) at (310:1.3cm) {};
\node (p6) at (350:1.3cm) {};
\node[rectangle,draw=black!20!white,fill=black!20!white](l1) at (50:1.05cm){$\gamma^2(v_1)$};
\node[rectangle,draw=black!20!white,fill=black!20!white](l2) at (130:1.05cm){$\gamma^2(v_2)$};
\node[rectangle,draw=black!20!white,fill=black!20!white](l3) at (180:0.95cm){$v_1$};
\node[rectangle,draw=black!20!white,fill=black!20!white](l4) at (240:0.95cm){$v_2$};
\filldraw[fill=black!50!white, draw=black, thick]
    (2.75,0) -- (3.35,0) -- (3.35,-0.1) -- (3.55,0.05) -- (3.35,0.2) -- (3.35,0.1) -- (2.75,0.1) -- (2.75,0);
\end{tikzpicture}
\hspace{0.5cm}
\begin{tikzpicture}[very thick,scale=1]
\tikzstyle{every node}=[circle, draw=black, fill=white, inner sep=0pt, minimum width=5pt];
\filldraw[fill=black!20!white, draw=black, thin, dashed](0,0)circle(1.6cm);
\node[rectangle,draw=black!20!white,fill=black!20!white](l5) at (290:0.81cm){$\gamma(v_1)$};
\node[rectangle,draw=black!20!white,fill=black!20!white](l6) at (360:0.71cm){$\gamma(v_2)$};
\node (p1) at (70:1.3cm) {};
\node (p2) at (110:1.3cm) {};
\node (p3) at (190:1.3cm) {};
\node (p4) at (230:1.3cm) {};
\node (p5) at (310:1.3cm) {};
\node (p6) at (350:1.3cm) {};
\node (p7) at (90:2.2cm) {};
\node (p8) at (210:2.2cm) {};
\node (p9) at (330:2.2cm) {};
\node[rectangle,draw=black!20!white,fill=black!20!white](l1) at (50:1.05cm){$\gamma^2(v_1)$};
\node[rectangle,draw=black!20!white,fill=black!20!white](l2) at (130:1.05cm){$\gamma^2(v_2)$};
\node[rectangle,draw=black!20!white,fill=black!20!white](l3) at (180:0.95cm){$v_1$};
\node[rectangle,draw=black!20!white,fill=black!20!white](l4) at (240:0.95cm){$v_2$};
\node[rectangle,draw=white,fill=white](l7) at (100:2.3cm){$z$};
\node[rectangle,draw=white,fill=white](l8) at (210:2.5cm){$v$};
\node[rectangle,draw=white,fill=white](l9) at (330:2.5cm){$w$};
\draw[thin](p7)--(p1);
\draw[ultra thick](p7)--(p2);
\draw[ultra thick](p8)--(p3);
\draw[dashed](p8)--(p4);
\draw[dashed](p9)--(p5);
\draw[thin](p9)--(p6);
\end{tikzpicture}
\end{center}
\vspace{-0.3cm}
\caption{\emph{Construction of a $(\mathcal{C}_{3},\Phi)$ 3Tree2 partition of $G$ in the case where $G$ is a $(\mathcal{C}_{3},\Phi_{k-1})$ vertex addition of $G_{k-1}$.}}
\label{C3treevert}
\end{figure}
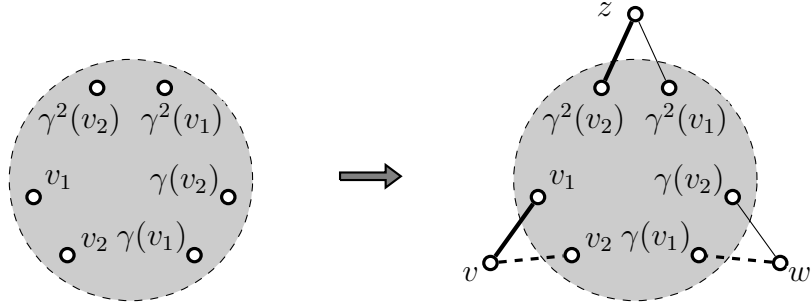

Suppose next that $G$ is a $(\mathcal{C}_{3},\Phi_{k-1})$ edge split on $(\{v_{1},v_{2}\},\{\gamma(v_{1}),\gamma(v_{2})\},\{\gamma^2(v_{1}),\gamma^2(v_{2})\});(v,w,z)$ of $G_{k-1}$ with $E(G_{k})=\big(E(G_{k-1})\setminus\big\{\{v_{1},v_{2}\},\{\gamma(v_{1}),\gamma(v_{2})\},\{\gamma^2(v_{1}),
\gamma^2(v_{2})\}\big\}\big)\cup \big\{\{v,v_{i}\}|\, i=1,2,3\big\}\cup \big\{\{w,\gamma(v_{i})\}|\, i=1,2,3\big\}\cup \big\{\{z,\gamma^2(v_{i})\}|\, i=1,2,3\big\}$, where $w=\gamma(v)$ and $z=\gamma^2(v)$. Wlog we may assume that $\{v_{1},v_{2}\}\in E\big(T^{(k-1)}_{0}\big)$. Then $\{\gamma(v_{1}),\gamma(v_{2})\}\in E\big(T^{(k-1)}_{1}\big)$ and $\{\gamma^2(v_{1}),\gamma^2(v_{2})\}\in E\big(T^{(k-1)}_{2}\big)$. Note that $v_{3}$ belongs to a tree $T^{(k-1)}_{l}$, where $l\neq 0$. Suppose $v_{3}\in T^{(k-1)}_{1}$. Then $\gamma(v_{3})\in T^{(k-1)}_{2}$ and $\gamma^2(v_{3})\in T^{(k-1)}_{0}$.
\begin{figure}[htp]
\begin{center}
\begin{tikzpicture}[very thick,scale=1]
\tikzstyle{every node}=[circle, draw=black, fill=white, inner sep=0pt, minimum width=5pt];
\filldraw[fill=black!20!white, draw=black, thin, dashed](0,0)circle(1.8cm);
\node (p1) at (70:1.5cm) {};
\node (p2) at (110:1.5cm) {};
\node (p3) at (190:1.5cm) {};
\node (p4) at (230:1.5cm) {};
\node (p5) at (310:1.5cm) {};
\node (p6) at (350:1.5cm) {};
\node (p7) at (90:0.7cm) {};
\node (p8) at (210:0.7cm) {};
\node (p9) at (330:0.7cm) {};
\node[rectangle,draw=black!20!white,fill=black!20!white](l1) at (50:1.25cm){$\gamma^2(v_1)$};
\node[rectangle,draw=black!20!white,fill=black!20!white](l2) at (130:1.25cm){$\gamma^2(v_2)$};
\node[rectangle,draw=black!20!white,fill=black!20!white](l3) at (174:1.35cm){$v_1$};
\node[rectangle,draw=black!20!white,fill=black!20!white](l4) at (246:1.35cm){$v_2$};
\node[rectangle,draw=black!20!white,fill=black!20!white](l5) at (287:1.25cm){$\gamma(v_1)$};
\node[rectangle,draw=black!20!white,fill=black!20!white](l6) at (368:1.3cm){$\gamma(v_2)$};
\node[rectangle,draw=black!20!white,fill=black!20!white](l7) at (90:0.3cm){$\gamma^2(v_3)$};
\node[rectangle,draw=black!20!white,fill=black!20!white](l8) at (235:0.9cm){$v_3$};
\node[rectangle,draw=black!20!white,fill=black!20!white](l9) at (305:0.88cm){$\gamma(v_3)$};
\draw[thin](p1)--(p2);
\draw[ultra thick](p3)--(p4);
\draw[dashed](p5)--(p6);
\filldraw[fill=black!50!white, draw=black, thick]
    (2.75,0) -- (3.35,0) -- (3.35,-0.1) -- (3.55,0.05) -- (3.35,0.2) -- (3.35,0.1) -- (2.75,0.1) -- (2.75,0);
\end{tikzpicture}
\hspace{0.5cm}
\begin{tikzpicture}[very thick,scale=1]
\tikzstyle{every node}=[circle, draw=black, fill=white, inner sep=0pt, minimum width=5pt];
\filldraw[fill=black!20!white, draw=black, thin, dashed](0,0)circle(1.8cm);
\node (p1) at (70:1.5cm) {};
\node (p2) at (110:1.5cm) {};
\node (p3) at (190:1.5cm) {};
\node (p4) at (230:1.5cm) {};
\node (p5) at (310:1.5cm) {};
\node (p6) at (350:1.5cm) {};
\node (a7) at (90:2.4cm) {};
\node (a8) at (210:2.4cm) {};
\node (a9) at (330:2.4cm) {};
\node (p7) at (90:0.7cm) {};
\node (p8) at (210:0.7cm) {};
\node (p9) at (330:0.7cm) {};
\node[rectangle,draw=black!20!white,fill=black!20!white](l1) at (50:1.25cm){$\gamma^2(v_1)$};
\node[rectangle,draw=black!20!white,fill=black!20!white](l2) at (130:1.25cm){$\gamma^2(v_2)$};
\node[rectangle,draw=black!20!white,fill=black!20!white](l3) at (174:1.35cm){$v_1$};
\node[rectangle,draw=black!20!white,fill=black!20!white](l4) at (246:1.35cm){$v_2$};
\node[rectangle,draw=black!20!white,fill=black!20!white](l5) at (287:1.25cm){$\gamma(v_1)$};
\node[rectangle,draw=black!20!white,fill=black!20!white](l6) at (368:1.3cm){$\gamma(v_2)$};
\node[rectangle,draw=black!20!white,fill=black!20!white](l7) at (90:0.3cm){$\gamma^2(v_3)$};
\node[rectangle,draw=black!20!white,fill=black!20!white](l8) at (235:0.9cm){$v_3$};
\node[rectangle,draw=black!20!white,fill=black!20!white](l9) at (305:0.88cm){$\gamma(v_3)$};
\node[rectangle,draw=white,fill=white](l7) at (97:2.5cm){$z$};
\node[rectangle,draw=white,fill=white](l8) at (210:2.7cm){$v$};
\node[rectangle,draw=white,fill=white](l9) at (330:2.7cm){$w$};
\draw[thin](a7)--(p1);
\draw[thin](a7)--(p2);
\draw[ultra thick](a7)--(p7);
\draw[ultra thick](a8)--(p3);
\draw[ultra thick](a8)--(p4);
\draw[dashed](a8)--(p8);
\draw[dashed](a9)--(p5);
\draw[dashed](a9)--(p6);
\draw[thin](a9)--(p9);
\end{tikzpicture}
\end{center}
\vspace{-0.3cm}
\caption{\emph{Construction of a $(\mathcal{C}_{3},\Phi)$ 3Tree2 partition of $G$ in the case where $G$ is a $(\mathcal{C}_{3},\Phi_{k-1})$ edge split of $G_{k-1}$.}}
\label{C3treeedge}
\end{figure}
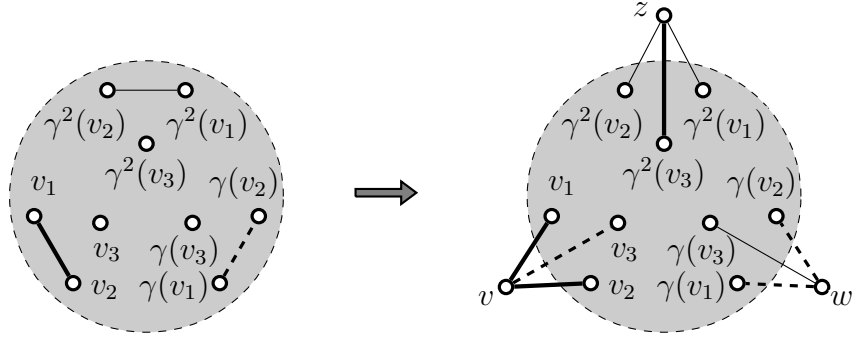
So, if we define $T^{(k)}_{0}$ to be the tree with
\begin{eqnarray}V\big(T^{(k)}_{0}\big)&=&V\big(T^{(k-1)}_{0}\big)\cup\{v,z\}\nonumber\\ E\big(T^{(k)}_{0}\big)&=&\big(E\big(T^{(k-1)}_{0}\big)\setminus \{v_{1},v_{2}\}\big)\cup\big\{\{v,v_{1}\},\{v,v_{2}\},\{z,\gamma^2(v_{3})\}\big\}\textrm{,}\nonumber\end{eqnarray}
$T^{(k)}_{1}$ to be the tree with
\begin{eqnarray}V\big(T^{(k)}_{1}\big)&=&V\big(T^{(k-1)}_{1}\big)\cup\{w,v\}\nonumber\\ E\big(T^{(k)}_{1}\big)&=&\big(E\big(T^{(k-1)}_{1}\big)\setminus \{\gamma(v_{1}),\gamma(v_{2})\}\big)\nonumber\\&&\cup\big\{\{w,\gamma(v_{1})\},\{w,\gamma(v_{2})\},\{v,v_{3}\}\big\}\textrm{,}
\nonumber\end{eqnarray}
and $T^{(k)}_{2}$ to be the tree with
\begin{eqnarray}V\big(T^{(k)}_{2}\big)&=&V\big(T^{(k-1)}_{2}\big)\cup\{z,w\}\nonumber\\ E\big(T^{(k)}_{2}\big)&=&\big(E\big(T^{(k-1)}_{2}\big)\setminus \{\gamma^2(v_{1}),\gamma^2(v_{2})\}\big)\nonumber\\&&\cup\big\{\{z,\gamma^2(v_{1})\},\{z,\gamma^2(v_{2})\},\{w,\gamma(v_{3})\}\big\}\textrm{,}
\nonumber\end{eqnarray}
then $\big\{E\big(T^{(k)}_{0}\big),E\big(T^{(k)}_{1}\big),E\big(T^{(k)}_{2}\big)\big\}$ is a $(\mathcal{C}_{3},\Phi)$ 3Tree2 partition of $G$. If $v_{3}\in T^{(k-1)}_{2}$, then we obtain a $(\mathcal{C}_{3},\Phi)$ 3Tree2 partition of $G$ in an analogous manner.

Finally, suppose that $G$ is a $(\mathcal{C}_{3},\Phi_{k-1})$ $\Delta$ extension by $(v\,w\,z)$ of $G_{k-1}$ with $E(G)=E(G_{k-1})\cup \big\{\{v,w\},\{w,z\},\{z,v\},\{v,v_{0}\},\{w,\gamma(v_{0})\},\{z,\gamma^2(v_{0})\}\big\}$, where $w=\gamma(v)$ and $z=\gamma^2(v)$. Wlog we may assume that $v_{0}\in V\big(T^{(k-1)}_{0}\big)$. Then $\gamma(v_{0})\in V\big(T^{(k-1)}_{1}\big)$ and $\gamma^2(v_{0})\in V\big(T^{(k-1)}_{2}\big)$.
\begin{figure}[htp]
\begin{center}
\begin{tikzpicture}[very thick,scale=1]
\tikzstyle{every node}=[circle, draw=black, fill=white, inner sep=0pt, minimum width=5pt];
\filldraw[fill=black!20!white, draw=black, thin, dashed](0,0)circle(1.3cm);
\node (p1) at (30:1cm) {};
\node (p2) at (150:1cm) {};
\node (p3) at (270:1cm) {};
\node[rectangle,draw=black!20!white,fill=black!20!white](l1) at (10:0.72cm){$\gamma(v_0)$};
\node[rectangle,draw=black!20!white,fill=black!20!white](l2) at (170:0.69cm){$\gamma^2(v_0)$};
\node[rectangle,draw=black!20!white,fill=black!20!white](l3) at (270:0.7cm){$v_0$};
\filldraw[fill=black!50!white, draw=black, thick]
    (2.75,0) -- (3.35,0) -- (3.35,-0.1) -- (3.55,0.05) -- (3.35,0.2) -- (3.35,0.1) -- (2.75,0.1) -- (2.75,0);
\end{tikzpicture}
\hspace{0.5cm}
\begin{tikzpicture}[very thick,scale=1]
\tikzstyle{every node}=[circle, draw=black, fill=white, inner sep=0pt, minimum width=5pt];
\filldraw[fill=black!20!white, draw=black, thin, dashed](0,0)circle(1.3cm);
\node (p1) at (30:1cm) {};
\node (p2) at (150:1cm) {};
\node (p3) at (270:1cm) {};
\node (p4) at (90:2.8cm) {};
\node (p5) at (210:2.8cm) {};
\node (p6) at (330:2.8cm) {};
\node[rectangle,draw=black!20!white,fill=black!20!white](l1) at (10:0.7cm){$\gamma(v_0)$};
\node[rectangle,draw=black!20!white,fill=black!20!white](l2) at (170:0.69cm){$\gamma^2(v_0)$};
\node[rectangle,draw=black!20!white,fill=black!20!white](l3) at (270:0.7cm){$v_0$};
\node[rectangle,draw=white,fill=white](l4) at (97:2.8cm){$z$};
\node[rectangle,draw=white,fill=white](l5) at (208:3.11cm){$v$};
\node[rectangle,draw=white,fill=white](l6) at (332:3.11cm){$w$};
\draw[thin](p4)--(p5);
\draw[ultra thick](p5)--(p6);
\draw[dashed](p6)--(p4);
\draw[thin](p4)--(p2);
\draw[ultra thick](p5)--(p3);
\draw[dashed](p6)--(p1);
\end{tikzpicture}
\end{center}
\vspace{-0.3cm}
\caption{\emph{Construction of a $(\mathcal{C}_{3},\Phi)$ 3Tree2 partition of $G$ in the case where $G$ is a $(\mathcal{C}_{3},\Phi_{k-1})$ $\Delta$ extension of $G_{k-1}$.}}
\label{C3treedelta}
\end{figure}
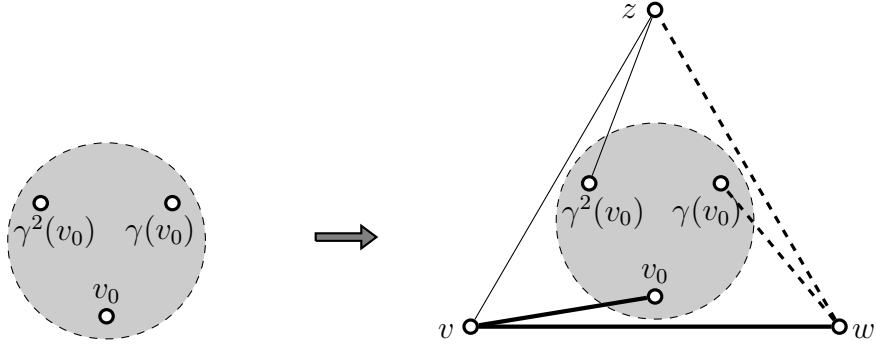
So, if we define $T^{(k)}_{0}$ to be the tree with \begin{eqnarray}V\big(T^{(k)}_{0}\big)&=&V\big(T^{(k-1)}_{0}\big)\cup\{v,w\}\nonumber\\ E\big(T^{(k)}_{0}\big)&=&E\big(T^{(k-1)}_{0}\big)\cup\big\{\{v,v_{0}\},\{v,w\}\big\}\textrm{,}\nonumber\end{eqnarray} $T^{(k)}_{1}$ to be the tree with
\begin{eqnarray}V\big(T^{(k)}_{1}\big)&=&V\big(T^{(k-1)}_{1}\big)\cup\{w,z\}\nonumber\\ E\big(T^{(k)}_{1}\big)&=&E\big(T^{(k-1)}_{1}\big)\cup\big\{\{w,\gamma(v_{0})\},\{w,z\}\big\}\textrm{,}\nonumber\end{eqnarray} and $T^{(k)}_{2}$ to be the tree with
\begin{eqnarray}V\big(T^{(k)}_{2}\big)&=&V\big(T^{(k-1)}_{2}\big)\cup\{v,z\}\nonumber\\ E\big(T^{(k)}_{2}\big)&=&E\big(T^{(k-1)}_{2}\big)\cup\big\{\{z,\gamma^2(v_{0})\},\{z,v\}\big\}\textrm{,}\nonumber\end{eqnarray} then $\big\{E\big(T^{(k)}_{0}\big),E\big(T^{(k)}_{1}\big),E\big(T^{(k)}_{2}\big)\big\}$ is a $(\mathcal{C}_{3},\Phi)$ 3Tree2 partition of $G$. $\square$

In order to show that $(iv)$ implies $(i)$ in Theorem \ref{C3char} we use an approach that is in the style of Tay's proof (see \cite{tay}) of Crapo's original result. This requires the notion of a `frame', i.e., a generalized notion of a framework that allows joints to be located at the same point in space, even if their corresponding vertices are adjacent.

\begin{defin}
\label{frame}
\emph{\cite{tay} Let $G$ be a graph with $V(G)=\{v_{1},v_{2},\ldots,v_{n}\}$. A \emph{frame} in $\mathbb{R}^2$ is a triple $(G,p,q)$, where $p:V(G)\to \mathbb{R}^2$ and $q:E(G)\to \mathbb{R}^2 \setminus\{0\}$ are maps with the property that for all $\{v_{i},v_{j}\}\in E(G)$ there exists a scalar  $\lambda_{ij}\in \mathbb{R}$ (which is possibly zero) such that $p(v_{i})- p(v_{j})=\lambda_{ij}q(\{v_{i},v_{j}\})$.}
\end{defin}

\begin{defin}
\label{framerigmatrix}
\emph{The \emph{generalized rigidity matrix} of a frame $(G,p,q)$ in $\mathbb{R}^2$ is the matrix \begin{displaymath} \mathbf{R}(G,p,q)=\left(
\begin{array} {ccccccccccc }
& & & & & \vdots & & & & & \\
0 & \ldots & 0 &  q(\{v_{i},v_{j}\}) &0 &\ldots &0 &  -q(\{v_{i},v_{j}\}) &0 &\ldots &
0\\ & & & & & \vdots & & & & &\end{array}
\right)\textrm{,}\end{displaymath} i.e., for each edge $\{v_{i},v_{j}\}\in E(G)$, $\mathbf{R}(G,p,q)$ has the row with
$\big(q(\{v_{i},v_{j}\})\big)_{1}$ and $\big(q(\{v_{i},v_{j}\})\big)_{2}$ in the columns $2i-1$ and $2i$, $-\big(q(\{v_{i},v_{j}\})\big)_{1}$ and $-\big(q(\{v_{i},v_{j}\})\big)_{2}$ in
the columns $2(j-1)$ and $2j$, and $0$ elsewhere.\\\indent We say that $(G,p,q)$ is \emph{independent} if $\mathbf{R}(G,p,q)$ has linearly independent rows.}
\end{defin}

\begin{remark}
\label{indeprem}
\emph{If $(G,p,q)$ is a frame with the property that $p(v_{i})\neq p(v_{j})$ whenever $\{v_{i},v_{j}\}\in E(G)$, then we obtain the rigidity matrix of the framework $(G,p)$ by multiplying each row of $\mathbf{R}(G,p,q)$ by its corresponding scalar $\lambda_{ij}$. Therefore, if $(G,p,q)$ is independent, so is $(G,p)$.}
\end{remark}

\begin{lemma}\label{3ivtoi} Let $G$ be a graph with $|V(G)|\geq 3$, $\mathcal{C}_{3}=\{Id,C_{3},C_{3}^2\}$ be a symmetry group in dimension $2$, and $\Phi:\mathcal{C}_{3}\to \textrm{Aut}(G)$ be a homomorphism. If $G$ has a proper $(\mathcal{C}_{3},\Phi)$ 3Tree2 partition, then $\mathscr{R}_{(G,\mathcal{C}_{3},\Phi)}\neq \emptyset$ and $G$ is $(\mathcal{C}_{3},\Phi)$-generically isostatic.
\end{lemma}

\textbf{Proof.} Suppose $G$ has a proper $(\mathcal{C}_{3},\Phi)$ 3Tree2 partition $\{E(T_{0}),E(T_{1}),E(T_{2})\}$. By Theorem \ref{symgenrigthm}, it suffices to find some framework $(G,p)\in \mathscr{R}_{(G,\mathcal{C}_{3},\Phi)}$ that is isostatic. Since $G$ has a 3Tree2 partition, $G$ satisfies the count $|E(G)|=2|V(G)|-3$ (see Remark \ref{3t2remk}), and hence, by Theorem \ref{isoequiv}, it suffices to find a map $p:V(G)\to \mathbb{R}^2$ such that $(G,p)\in \mathscr{R}_{(G,\mathcal{C}_{3},\Phi)}$ is independent. In the following, we again denote $\Phi(C_{3})$ by $\gamma$ and $\Phi(C_{3}^2)$ by $\gamma^2$.\\\indent Let $e_{0}=(0,0)$, $e_{1}=(1,0)$, and $e_{2}=(\frac{1}{2},\frac{\sqrt{3}}{2})$. Also, for $i=0,1,2$, let $V_{i}$ be the set of vertices of $G$ that are not in $V(T_{i})$, and let $(G,p,q)$ be the frame with $p:V(G)\to \mathbb{R}^2$ and $q:E(G)\to \mathbb{R}^2$ defined by
\begin{eqnarray}
p(v) & = & e_{i} \quad \textrm{if } v\in V_{i} \nonumber\\
q(b) & = & \left\{ \begin{array}{lll}
e_{2}-e_{1} &  =(-\frac{1}{2},\frac{\sqrt{3}}{2}) & \textrm{if } b\in E(T_{0})\\
e_{0}-e_{2} &  =(-\frac{1}{2},-\frac{\sqrt{3}}{2}) & \textrm{if } b\in E(T_{1})\\
e_{1}-e_{0} &  =(1,0) & \textrm{if } b\in E(T_{2})
\end{array} \right.\textrm{.}\nonumber
\end{eqnarray}

\begin{figure}[htp]
\begin{center}
\begin{tikzpicture}[very thick,scale=1]
\tikzstyle{every node}=[circle, draw=black, fill=black!30!white, inner sep=0pt, minimum width=10pt];
\draw[dashed](-0.13,1.732)--node[rectangle, draw=white,fill=white,left =5pt] {$T_1$}(-1.13,0);
\draw[dashed](0.13,1.732)--(-0.87,0);
\draw[thin](-1,0.1)--(1,0.1);
\draw[thin](-1,-0.1)--node[rectangle, draw=white,fill=white,below =5pt] {$T_2$}(1,-0.1);
\draw[ultra thick](0.13,1.732)--node[rectangle, draw=white,fill=white,right =5pt] {$T_0$}(1.13,0);
\draw[ultra thick](-0.13,1.732)--(0.87,0);
\node (p1) at (0,1.732) {};
\node (p2) at (-1,0) {};
\node (p3) at (1,0) {};
\node[rectangle,draw=white,fill=white](l1) at (0,2.2){$V_2$};
\node[rectangle,draw=white,fill=white](l2) at (-1.45,0){$V_0$};
\node[rectangle,draw=white,fill=white](l3) at (1.5,0){$V_1$};
\node[rectangle,draw=white,fill=white](l1) at (-0.5,1.8){$e_2$};
\node[rectangle,draw=white,fill=white](l2) at (-1,-0.4){$e_0$};
\node[rectangle,draw=white,fill=white](l3) at (1,-0.4){$e_1$};
\draw[dashed](p1)--(p2);
\draw[thin](p3)--(p2);
\draw[ultra thick](p1)--(p3);
\end{tikzpicture}
\end{center}
\vspace{-0.3cm}
\caption{\emph{The frame $(G,p,q)$.}}
\label{frameC3}
\end{figure}
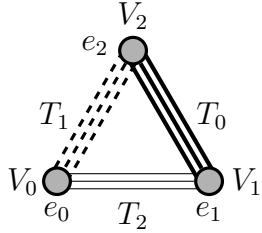

We claim that the generalized rigidity matrix $\mathbf{R}(G,p,q)$ has linearly independent rows. To see this, we first rearrange the columns of $\mathbf{R}(G,p,q)$ in such a way that we obtain the matrix $\mathbf{R'}(G,p,q)$ which has the $(2i-1)^{st}$ column of $\mathbf{R}(G,p,q)$ in its $i^{th}$ column and the $(2i)^{th}$ column of $\mathbf{R}(G,p,q)$ in its $(|V(G)|+i)^{th}$ column for $i=1,2,\ldots,|V(G)|$. Let $F_{b}$ denote the row vector of $\mathbf{R'}(G,p,q)$ that corresponds to the edge $b\in E(G)$. We then rearrange the rows of $\mathbf{R'}(G,p,q)$ in such a way that we obtain the matrix $\mathbf{R''}(G,p,q)$ which has the vectors $F_{b}$ with $b\in E(T_{0})$ in the rows $1,2,\ldots,|E(T_{0})|$, the vectors $F_{b}$ with $b\in E(T_{1})$ in the following $|E(T_{1})|$ rows, and the vectors $F_{b}$ with $b\in E(T_{2})$ in the last $|E(T_{2})|$ rows. So $\mathbf{R''}(G,p,q)$ is a matrix of the form
\begin{displaymath}\addtolength{\arraycolsep}{-0.8mm}
\left( \begin{array} {rrrrrrr|rrrrrrr}
 -\frac{1}{2}& & & \frac{1}{2} & & & & \frac{\sqrt{3}}{2} & & & -\frac{\sqrt{3}}{2} & & & \\
 & & & \vdots & & & &              & & & \vdots & & \\
 & & -\frac{1}{2} & & & & \frac{1}{2} &            & &  \frac{\sqrt{3}}{2}& & & & -\frac{\sqrt{3}}{2}\\
\hline
 & -\frac{1}{2} & & & \frac{1}{2} & & &  &  -\frac{\sqrt{3}}{2}& & &  \frac{\sqrt{3}}{2}&  &\\
 & & & \vdots & & & &              & & & \vdots & & & \\
 & & & -\frac{1}{2}& & \frac{1}{2} & & & & &  -\frac{\sqrt{3}}{2} & &  \frac{\sqrt{3}}{2}& \\
\hline
 1& & & & -1 & & &   & & & & & & \\
 & & & \vdots & & & &       & & & \mathbf{0}& & & \\
 & & 1 & & & -1 & &     & & & & & &
\end{array}
\right)\textrm{.}
\end{displaymath}

Clearly, $\mathbf{R}(G,p,q)$ has a row dependency if and only if $\mathbf{R''}(G,p,q)$ does. Suppose $\mathbf{R''}(G,p,q)$ has a row dependency of the form \begin{displaymath}\sum_{b\in E(G)}\alpha_{b}F_{b}=0\textrm{,}\end{displaymath} where $\alpha_{b}\neq 0$ for some $b\in E(T_{2})$. Then, since $T_{2}$ is a tree, we have \begin{displaymath}\sum_{b\in E(T_{2})}\alpha_{b}F_{b}\neq 0\textrm{.}\end{displaymath} Thus, there exists a vertex $v_{s}\in V(T_{2})$, $s\in\{1,2,\ldots,|V(G)|\}$, such that \begin{displaymath}\sum_{b\in E(T_{2})}\alpha_{b}(F_{b})_{s}=C\neq 0\textrm{.}\end{displaymath} Since $v_{s}\in V(T_{2})$, $v_{s}$ belongs to either $T_{0}$ or $T_{1}$, say wlog $v_{s}\in V(T_{1})$ and $v_{s}\notin V(T_{0})$. Therefore, $(F_{b})_{s}=0$ and $(F_{b})_{|V(G)|+s}=0$ for all $b\in E(T_{0})$ and \begin{displaymath}\sum_{b\in E(T_{1})}\alpha_{b}(F_{b})_{s}=-C\textrm{.}\end{displaymath} This says that \begin{displaymath}\sum_{b\in E(T_{1})}\alpha_{b}(F_{b})_{|V(G)|+s}=\sum_{b\in E(G)}\alpha_{b}(F_{b})_{|V(G)|+s}=-\sqrt{3}C\neq 0\textrm{,}\end{displaymath} a contradiction. Therefore, if $\sum_{b\in E(G)}\alpha_{b}F_{b}=0$ is a row dependency of $\mathbf{R''}(G,p,q)$, then $\alpha_{b}= 0$ for all $b\in E(T_{2})$.\\\indent So, it is now only left to show that the matrix $\mathbf{\widetilde{R}}(G,p,q)$ which is obtained from $\mathbf{R''}(G,p,q)$ by deleting those rows of $\mathbf{R''}(G,p,q)$ that correspond to the edges of $T_{2}$ has linearly independent rows. This can be done by multiplying $\mathbf{\widetilde{R}}(G,p,q)$ by  appropriate matrices of basis transformations and then using arguments analogous to those above. So, as claimed, the frame $(G,p,q)$ is independent.
\\\indent Now, if $(G,p)$ is not a framework, then we need to symmetrically pull apart those joints of $(G,p,q)$ that have the same location $e_{i}$ in $\mathbb{R}^2$ and whose vertices are adjacent. So, wlog suppose $|V_{0}|\geq 2$. Then we also have $|V_{0}|=|V_{1}|=|V_{2}|\geq 2$, because $\{E(T_{0}),E(T_{1}),E(T_{2})\}$ is a $(\mathcal{C}_{3},\Phi)$ 3Tree2 partition of $G$. Since $\{E(T_{0}),E(T_{1}),E(T_{2})\}$ is proper, one of $\langle V_{0} \rangle\cap T_{i}$, $i=1,2$, say wlog $\langle V_{0} \rangle\cap T_{2}$, is not connected, and hence $\langle V_{1} \rangle\cap T_{0}$ and $\langle V_{2} \rangle\cap T_{1}$ are also not connected. Let $A$ be the set of vertices in one of the components of $\langle V_{0} \rangle\cap T_{2}$ and $\gamma(A)$ and $\gamma^2(A)$ be the vertex sets of the corresponding components of $\langle V_{1} \rangle\cap T_{0}$ and $\langle V_{2} \rangle\cap T_{1}$, respectively. For $t\in \mathbb{R}$, we define $p_{t}:V(G)\to \mathbb{R}^2$ and $q_{t}:E(G)\to \mathbb{R}^2$ by
\begin{eqnarray}
p_{t}(v)& = & \left\{ \begin{array}{ll}
(-\frac{1}{2}t,-\frac{\sqrt{3}}{2}t) & \textrm{if } v\in A\\
(1+t,0) & \textrm{if } v\in \gamma(A)\\
\big(\frac{1}{2}(1-t),\frac{\sqrt{3}}{2}(1+t)\big) & \textrm{if } v\in \gamma^2(A)\\
p(v) & \textrm{otherwise}
\end{array}\right. \nonumber\\
q_{t}(b)& = & \left\{ \begin{array}{ll}
(1+\frac{1}{2}t,\frac{\sqrt{3}}{2}t) & \textrm{if } b\in E_{A,V_{1}\setminus{\gamma(A)}}\\
(1+\frac{3}{2}t,\frac{\sqrt{3}}{2}t) & \textrm{if } b\in E_{A,\gamma(A)}\\
(-\frac{1}{2}-t,\frac{\sqrt{3}}{2}) & \textrm{if } b\in E_{\gamma(A),V_{2}\setminus{\gamma^2(A)}}\\
\big(-\frac{1}{2}-\frac{3}{2}t,\frac{\sqrt{3}}{2}(1+t)\big) & \textrm{if } b\in E_{\gamma(A),\gamma^2(A)}\\
\big(-\frac{1}{2}(1-t),-\frac{\sqrt{3}}{2}(1+t)\big) & \textrm{if } b\in E_{\gamma^2(A),V_{0}\setminus{A}}\\
(-\frac{1}{2},-\frac{\sqrt{3}}{2}-\sqrt{3}t) & \textrm{if } b\in E_{\gamma^2(A),A}\\
q(b) & \textrm{otherwise}
\end{array} \right.\textrm{,}\nonumber
\end{eqnarray}
where for disjoint sets $X,Y\in V(G)$, $E_{X,Y}$ denotes the set of edges of $G$ incident with a vertex in $X$ and a vertex in $Y$.
\begin{figure}[htp]
\begin{center}
\begin{tikzpicture}[very thick,scale=1]
\tikzstyle{every node}=[circle, draw=black, fill=black!30!white, inner sep=0pt, minimum width=10pt];
\draw[dashed](-0.13,1.732)--(-1.13,0);
\draw[dashed](0.13,1.732)--(-0.87,0);
\draw[thin](-1,0.1)--(1,0.1);
\draw[thin](-1,-0.1)--(1,-0.1);
\draw[ultra thick](0.13,1.732)--(1.13,0);
\draw[ultra thick](-0.13,1.732)--(0.87,0);
\node (p1) at (0,1.732) {};
\node (p2) at (-1,0) {};
\node (p3) at (1,0) {};
\node[rectangle,draw=white,fill=white](l1) at (1.32,1.9){$V_2\setminus \gamma^2(A)$};
\node[rectangle,draw=white,fill=white](l2) at (-2,0){$V_0\setminus A$};
\node[rectangle,draw=white,fill=white](l3) at (1,-0.8){$V_1\setminus \gamma(A)$};
\node[rectangle,draw=white,fill=white](l1) at (-0.45,1.7){$e_2$};
\node[rectangle,draw=white,fill=white](l2) at (-0.8,-0.4){$e_0$};
\node[rectangle,draw=white,fill=white](l3) at (1.3,0.3){$e_1$};
\node[rectangle,draw=white,fill=white](l1) at (-1.9,-0.866){$A$};
\node[rectangle,draw=white,fill=white](l2) at (2.67,0){$\gamma(A)$};
\node[rectangle,draw=white,fill=white](l3) at (-1.25,2.6){$\gamma^2(A)$};
\draw[dashed](p1)--(p2);
\draw[thin](p3)--(p2);
\draw[ultra thick](p1)--(p3);
\draw[ultra thick](-0.33,2.63)--node[rectangle, draw=white,fill=white,right =6pt] {$T_0$}(2,0);
\draw[dashed](-0.5,2.598)--node[rectangle, draw=white,fill=white,left =5pt] {$T_1$}(-1.65,-0.866);
\draw[thin](-1.5,-0.866)--node[rectangle, draw=white,fill=white,below left=5pt] {$T_2$}(2,-0.15);
\draw[ultra thick](0.1,1.732)--(2,0);
\draw[thin](1,-0.1)--(-1.5,-0.866);
\draw[dashed](-0.5,2.598)--(-1.1,0);
\draw[ultra thick](-0.4,2.598)--(0.1,1.732);
\draw[ultra thick](-0.6,2.598)--(-0.1,1.732);
\draw[thin](0.9,0.08)--(2,0.08);
\draw[thin](0.9,-0.08)--(2,-0.08);
\draw[dashed](-1.6,-0.9)--(-1.1,0);
\draw[dashed](-1.4,-0.9)--(-0.9,0);
\node (p4) at (-0.5,2.598) {};
\node (p5) at (-1.5,-0.866) {};
\node (p6) at (2,0) {};
\node (p1) at (0,1.732) {};
\node (p2) at (-1,0) {};
\node (p3) at (1,0) {};
\end{tikzpicture}
\end{center}
\vspace{-0.3cm}
\caption{\emph{The frame $(G,p_{t},q_{t})$.}}
\label{frameC3ext}
\end{figure}
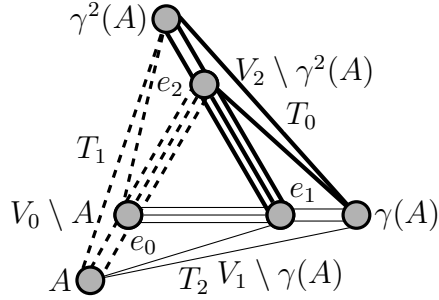
Then $(G,p_{t},q_{t})=(G,p,q)$ if $t=0$. Now, if we let $t'$ be an indeterminate, then the rows of $(G,p_{t'},q_{t'})$ are linearly dependent (over the quotient field of $\mathbb{R}[t']$) if and only if the determinants of all $|E(G)|\times |E(G)|$ submatrices of $(G,p_{t'},q_{t'})$ are identically zero. These determinants are polynomials in $t'$. Thus, the set of all $t\in \mathbb{R}$ with the property that $\mathbf{R}(G,p_{t},q_{t})$ has a non-trivial row dependency is a variety $F$ whose complement, if non-empty, is a dense open set. Since $t=0$ is in the complement of $F$ we can conclude that for almost all $t$, $(G,p_{t},q_{t})$ is independent. Therefore, there exists a $t_{0}\in\mathbb{R}$, $t_{0}\neq 0$, such that the frame $(G,p_{t_{0}},q_{t_{0}})$ is independent. This process can be continued until we obtain an independent frame $(G,\hat{p},\hat{q})$ with $\hat{p}(u)\neq \hat{p}(v)$ for all $\{u,v\}\in E(G)$. Then, by Remark \ref{indeprem}, $(G,\hat{p})$ is an independent framework and the right translation of $(G,\hat{p})$ yields an independent framework in the set $\mathscr{R}_{(G,\mathcal{C}_{3},\Phi)}$. $\square$

Lemmas \ref{3itoii}, \ref{3iitoiii}, \ref{3iiitoiv}, and \ref{3ivtoi} provide a complete proof for Theorem \ref{C3char}.

As shown in \cite{BS4}, there also exists a direct \emph{geometric} proof for the fact that condition $(iii)$ implies condition $(i)$ in Theorem \ref{C3char}, i.e., that the existence of a $(\mathcal{C}_{3},\Phi)$ construction sequence for $G$ implies that $\mathscr{R}_{(G,\mathcal{C}_{3},\Phi)}\neq \emptyset$ and that $G$ is $(\mathcal{C}_{3},\Phi)$-generically isostatic. By generalizing the basic geometric techniques used in this proof, we can construct classes of $(S,\Phi)$-generically isostatic graphs for a variety of symmetry groups $S$. These techniques also allow us to prove (or at least conjecture) characterizations of $(S,\Phi)$-generically isostatic graphs in situations where symmetric tree partitions are too complex. Moreover, they provide significant results for $(S,\Phi)$-generically independent graphs.

An immediate consequence of the symmetrized Laman's Theorems for $\mathcal{C}_{3}$, $\mathcal{C}_{2}$, and $\mathcal{C}_{s}$ (and the analogous conjectures for $\mathcal{C}_{2v}$ and $\mathcal{C}_{3v}$) is that there is (would be) a polynomial time algorithm to determine whether a given graph $G$ is $(S,\Phi)$-generically isostatic. In fact, although the Laman conditions involve an exponential number of subgraphs of $G$, there are several algorithms that determine whether they hold in $c|V(G)||E(G)|$ steps, where $c$ is a constant. The pebble game (\cite{henjac}) is an example for such an algorithm. The additional symmetry conditions for the number of fixed structural components can trivially be checked in constant time, from the graph automorphisms.

\section{Further work}

 \subsection{Pinned frameworks}

In mechanical and structural engineering, one is often interested in the rigidity and flexibility properties of \emph{pinned} frameworks, i.e., frameworks that have some of their joints firmly anchored (`pinned') to the ground (see, for example, \cite{FGsymmax, Assur1, Assur2}). Using the techniques presented in \cite{cfgsw,FGsymmax,  BS4}, it is straightforward to show that an isostatic  symmetric \emph{pinned} framework $(G,p)$ must again satisfy some very simply stated restrictions on the number of (unpinned) joints and bars  of $(G,p)$ that are fixed by various symmetry operations of $(G,p)$. While there are only six possible point groups that allow isostatic frameworks in the plane, it turns out that an isostatic \emph{pinned} framework can be constructed for \emph{any} point group in dimension 2.\\\indent We conjecture that the standard Laman-type conditions for a pinned graph $G$ (see \cite{Assur1}, for example), together with the additional necessary conditions concerning the number of fixed   joints and bars, are also sufficient for pinned 2-dimensional realizations of $G$ which are as generic as possible subject to the given symmetry constraints to be isostatic.\\\indent In particular, for the symmetry groups $\mathcal{C}_2$, $\mathcal{C}_3$, and $\mathcal{C}_s$ in dimension 2, we claim that the techniques of this paper extend directly to proofs of the corresponding symmetric versions of Laman's Theorem for pinned frameworks.

\subsection{Frameworks in dimension $d>2$}

A combinatorial characterization of generically $d$-isostatic graphs in dimension $d> 2$ has not yet been found \cite{graver, gss, W2}. Recall from Section \ref{subsec:genrig}, however, that there are a number of inductive construction techniques which are known to preserve the generic rigidity properties of a graph (see also \cite{LMW, W1}).\\\indent
It is shown in \cite{BS4} that, unlike in dimension 2, symmetry in dimension $d>2$ induces extra conditions for a graph $G$ to be $(S,\Phi)$-generically isostatic beyond those of
\begin{itemize}
\item[(a)] $G$ being generically $d$-isostatic and
\item [(b)] the symmetry conditions derived in \cite{cfgsw} concerning the number of fixed structural components of $G$ (and of all symmetric subgraphs $H$ of $G$ with the full count $|E(H)|=d|V(H)|-\binom{d+1}{2}$).
\end{itemize}
We conjecture that `flatness' caused by symmetry is the only additional concern, and that it can be made into a finite set of added combinatorial conditions, for each symmetry group. See \cite {cfgsw, BS4, sww} for further details.

\subsection{Body-bar and body-hinge structures}

Faced with the difficulties of characterizing generically rigid graphs in dimension $d>2$, in contrast with the well developed theory in the plane, there has recently been a careful study of a special class of frameworks, the class of \emph{body-bar frameworks}. These structures have a basically complete combinatorial theory  which exhibits all the key theorems and algorithms of the well understood plane frameworks (see, for example, \cite{taybb, WWbb, W1}).\\\indent
For a body-bar framework $(G,p)$ that possesses non-trivial symmetries, joint work with S. Guest and W. Whiteley  shows that in addition to the conditions in Tay's Theorem (see \cite{taybb}), there exist further necessary conditions for $(G,p)$ to be isostatic \cite{gsw}. These can be formulated as restrictions on the number of bars and bodies that are fixed by various symmetry operations of $(G,p)$. While these extra conditions are analogous to the ones derived for bar and joint frameworks, the modified context holds the promise of converting them into necessary and sufficient conditions for an \emph{arbitrary-dimensional} body-bar realization of $G$ to be isostatic, provided that this realization is as generic as possible subject to the given symmetry constraints. These conjectures, as well as various additional conjectures concerning combinatorial characterizations of $d$-dimensional symmetric body-bar frameworks, are stated in \cite{gsw}.\\\indent
For the groups $\mathcal{C}_2$, $\mathcal{C}_3$, and $\mathcal{C}_s$ in dimension 2, these conjectures can readily be proven by  modeling a symmetric body-bar framework as a framework (in the sense of Definition \ref{framework}) with isostatic bar and joint bodies of required symmetry, and then applying the results of this paper.

An interesting special class of body-bar frameworks with some important applications in rigidity theory is the class of \emph{body-hinge frameworks} \cite{TW1, W1,W5, W2}. It is shown in \cite{TW1} that body-hinge realizations of a multigraph $G$ with generic hinge assignments are infinitesimally rigid if and only if body-bar realizations of $G$ with generic positions for the end-points of the bars are infinitesimally rigid. So, body-hinge frameworks have the same efficient algorithms for testing generic rigidity as body-bar frameworks \cite{W5}. Moreover, the Molecular Conjecture posed by T.-S. Tay and W. Whiteley in 1984 proposes that the even more special class of body-bar frameworks that arise in the models of molecular kinematics (i.e., the class of body-hinge frameworks that have all hinges of each body concurrent in a point) also have the same good combinatorial theory as general body-bar frameworks \cite{TW1, W5}, so that, under generic conditions, the efficient counting algorithms for body-bar frameworks also apply to molecular body-hinge frameworks.
\\\indent
Given certain symmetry constraints, we conjecture  that, analogously to the non-symmetric situation,  the results and conjectures in \cite{gsw} concerning symmetric-generic body-bar frameworks also translate directly to symmetric-generic body-hinge frameworks. We further conjecture that a symmetric version of the Molecular Conjecture holds, i.e., that under symmetric-generic conditions,  body-bar frameworks and molecular frameworks also possess the same rigidity properties.\\\indent We note that a number of biomolecules possess rotational symmetry, including a number of virus shells which exhibit the symmetry of the rotational icosahedral group. The potential for such applications, as well as for understanding human-built  structures which are designed to have symmetry, is a further motivation for giving explicit results for symmetric body-bar, body-hinge, and molecular structures.

\section*{Acknowledgements}

We would like to thank Walter Whiteley for numerous  interesting and helpful discussions.

\end{document}